\documentclass[11pt,a4paper]{article}
\usepackage{graphicx}
\usepackage{amsmath,amssymb}
\usepackage{epstopdf}
\usepackage[font=scriptsize]{caption}
\newtheorem{lemma}{\sc \bf Lemma}[section]
\newtheorem{propos}{\sc \bf Proposition}[section]
\newtheorem{theor}{\sc \bf Theorem}[section]
\newtheorem{corr}{\sc \bf Corollary}[section]
\newtheorem{remark}{\sc \bf Remark}[section]
\newtheorem{definition}{\sc \bf Definition}[section]
\newtheorem{example}{\sc \bf Example}[section]
\newtheorem{algorithm}{\sc \bf Algorithm}[section]
 \newcommand{\Rb}{\mathbb{R}}
\begin{document}

\thispagestyle{empty}

%
%
%
\title{A fuzzy method for solving fuzzy fractional differential equations based on
the generalized fuzzy Taylor expansion }

\author{ T. Allahviranloo$^{a,}$\footnote{Corresponding author,
E-mail: tofigh.allahviranloo@eng.bau.edu.tr;
allahviranloo@yahoo.com}, Z. Noeiaghdam$^b$ , S. Noeiaghdam$^{c,d}$,
S. Salahshour$^{a}$, Juan J. Nieto$^{e}$}
\date{}
\maketitle \vspace{-9mm}
\begin{center}
\scriptsize{ $^a$Faculty of Engineering and
Natural Sciences, Bahcesehir University, Istanbul, Turkey.\\
$^b$Department of Mathematics and Computer Science, Shahed University, Tehran, Iran. \\
$^c$South Ural State University, Lenin prospect 76, Chelyabinsk, 454080, Russian Federation.\\
$^d$Baikal School of BRICS, Irkutsk National Research Technical
University, Irkutsk, Russian Federation.\\
$^e$Instituto de Matem\'aticas, Departamento de Estat\'istica,
An\'alise Matem\'atica e Optimizaci\'on, Universidade de Santiago de
Compostela 15782, Santiago de Compostela, Spain.}
\end{center}
\date{}
%
%
%
\maketitle
\begin{abstract}
In many mathematical types of research, in order to solve the fuzzy fractional differential equations, we should transform these problems into crisp corresponding problems and by solving them the approximate solution can be obtained. The aim of this paper is to present a new direct method to solve the fuzzy fractional differential equations without this transformation. In this work, the fuzzy generalized Taylor expansion by using the sense of fuzzy Caputo fractional derivative for fuzzy-valued functions is presented. For solving fuzzy fractional differential equations, the fuzzy generalized Euler's method is applied. In order to show the accuracy and efficiency of the presented method, the local and global truncation errors are determined. Moreover, the consistency, the convergence and the stability of the generalized Euler's method are proved in detail. Eventually, the numerical examples, especially in the switching point case, show the flexibility and the capability of the presented method.

\vspace{.5cm}{\it Keywords:} Fuzzy fractional differential equations; Generalized fuzzy Taylor expansion; Generalized fuzzy Euler's method; Global truncation error; Local truncation error; Convergence; Stability.
\end{abstract}
\section{Introduction}
Fuzzy set theory is a powerful tool for modeling uncertain problems. Therefore, large varieties of natural phenomena have been modeled using fuzzy concepts. Particularly, the fuzzy fractional differential equation is a common model in different science, such as population models, evaluating weapon systems, civil engineering, and modeling electro-hydraulic. Hence the concept of the fractional derivative is a very important topic in fuzzy calculus. Therefore, the fuzzy fractional
differential equations have attracted lots of attention in
mathematics and engineering researches. First work
devoted to the subject of fuzzy fractional differential equations is
the paper by Agarwal et al. \cite{al}. They have defined the
Riemann-Liouville differentiability concept under the Hukuhara
differentiability to solve fuzzy fractional differential equations.

In recent years, fractional calculus has introduced as an applicable
topic to produce the accurate results of mathematical and
engineering problems such as aerodynamics and control systems,
signal processing, bio-mathematical problems and others \cite{al,
bl, ks, mi, zn1}.

Furthermore, fractional differential equations in the fuzzy case
\cite{al} have studied by many authors and they have solved by
various methods \cite{a0, ag2, fracdif3,fracdif4}. In
\cite{fracdif7} Hoa studied the fuzzy fractional
differential equations under Caputo gH-differentiability and in \cite{fracdif2}
Agarwal et al. had a survey on mentioned problem to show the its
relation with optimal control problems. Also, Long et al.
\cite{fracdif5} illustrated the solvability of fuzzy fractional
differential equations and Salahshour et al.\cite{fracdif6}
applied the fuzzy Laplace transforms to solve this problem.

There are many numerical methods to solve the fuzzy fractional
differential equations by transforming to crisp problems
\cite{NEW3,NEW4,NEW5}. In this paper, a new direct method is
introduced to solve the mentioned problem without changing to crisp
form. The Taylor expansion method is one of the famous and applicable
methods to solve the linear and non-linear problems \cite{nietonew,
man13, man10}. In this paper, the fuzzy generalized Taylor expansion based on the fuzzy Caputo fractional derivative is expanded. Then the Euler's method is applied to solve the fuzzy fractional differential equations. Also, the local and global truncation errors are considered and finally the consistency, the convergence and the stability of the generalized fuzzy Euler's method are demonstrated. Furthermore, some examples with the switching point are solved by using the presented method. The numerical results show the precision of the generalized Euler's method to solve the fuzzy fractional differential equations.
\section{Basic Concepts}
At first, the brief summary of the fuzzy details and some
preliminaries are revisited \cite{a1, a2, bs, go, mm, sb}.
\begin{definition}
Set $\mathbb{R}_{\mathcal{F}}=\{u:\mathbb{R}^n\rightarrow[0,1]$ such
that $u$ satisfies in the conditions I to IV $\}$
\begin{itemize}
\item [I.]$u$ is normal: there exists an $x_0\in \mathbb{R}^n$ such that $u(x_0)=1$,
\item [II.]$u$ is fuzzy convex: for $0\leq\lambda\leq 1,~u(\lambda_{x_1}+(1-\lambda_{x_2}))\geq\min\{u(x_1),u(x_2)\}$,
 \item[III.]$u$ is upper semi-continuous:  for any $x_0\in \mathbb{R}^n$, it holds that $u(x_0)\geq \lim_{x\rightarrow x_0^{\pm}}u(x)$,
 \item[IV.] $[u]_0=\overline{supp(u)}=cl\{x\in\mathbb{R}^n~\mid~u(x)>0\}$  is a compact
 subset,
\end{itemize}
is called the space of fuzzy numbers or the fuzzy numbers set. The
$r$-level set is $[u]_r=\{x\in\mathbb{R}^n \mid u(x)\geq r,~0<r\leq
1\}$. Then from I to IV, it follows that, the $r$-level sets of
$u\in\mathbb{R}_{\mathcal{F}}$ are nonempty, closed and bounded
intervals.
\end{definition}
\begin{definition}
A triangular fuzzy number is defined as a fuzzy set in
$\mathbb{R}_{\mathcal{F}}$, that is specified by an ordered triple
$u=(a,b,c)\in \mathbb{R}^3$ with $a\leq b\leq c$ such that
$\underline{u}(r)=a+(b-a)r$ (or lower bound of $u$) and
$\overline{u}(r)=c-(c-b)r$ (or upper bound of $u$) are the endpoints
of $r$-level sets for all $r\in[0,1]$.
\end{definition}
A crisp number $k$ is simply represented by
$\overline{u}(r)=\underline{u}(r)=k,~\ 0 \leq r \leq 1$ and called
singleton. For arbitrary $~u,~v\in \mathbb{R}_{\mathcal{F}}$ and
scalar $k,$ we might summarize the addition and the scalar
multiplication of two fuzzy numbers by

$ \triangleleft~~~~addition:~[u\oplus v]_r=[\underline{u}(r)+\underline{v}(r), \overline{u}(r)+\overline{v}(r)],$\\

$\triangleleft~~~~scalar~multiplication:~\left\{\begin{array}{ll}
[k\odot u]_r=[k\underline{u}(r),k\overline{u}(r)],~~~~k\geq 0, \\\\

[k\odot u]_r=[k\overline{u}(r),k\underline{u}(r)],~~~~k<0.
\end{array}
\right.$

 The Hausdorff distance between fuzzy numbers is given by
$\mathcal{H}:\mathbb{R}_{\mathcal{F}}\times
\mathbb{R}_{\mathcal{F}}\rightarrow \mathbb{R}^{+}\cup \{0\}$ as:
\[\mathcal{H}(u,v)=\sup_{0\leq
r\leq1}\max\{|\underline{u}(r)-\underline{v}(r)|,|\overline{u}(r)-\overline{v}(r)|\},\]
where $[u]_r= [\underline{u}(r),\overline{u}(r)],~[v]_r=
[\underline{v}(r),\overline{v}(r)]$. The metric space
$(\mathbb{R}_{\mathcal{F}} , \mathcal{H})$ is complete, separable
and locally compact  where the following conditions are valid for
metric $\mathcal{H}$:
\begin{itemize}
\item[I.]  $\mathcal{H}(u\oplus w, v\oplus w)=\mathcal{H}(u, v),~\forall u,~v,~w\in\mathbb{R}_{\mathcal{F}}.$
\item[II.]  $\mathcal{H}(\lambda u,\lambda v)=\mid \lambda\mid \mathcal{H}(u,v),~\forall\lambda\in\mathbb{R},~\forall u,~v\in \mathbb{R}_{\mathcal{F}}.$
\item[III.]  $\mathcal{H}(u\oplus v, w\oplus z)\leq \mathcal{H}(u,w)+\mathcal{H}(v,z),~\forall u,~v,~w,~z\in\mathbb{R}_{\mathcal{F}}.$
\end{itemize}
\begin{definition} Let$~u,v\in\mathbb{R}_{\mathcal{F}}$, if there exists$~w\in\mathbb{R}_{\mathcal{F}},~$such
that$~u=v+w,~$then$~w~$is called the Hukuhara difference
(H-difference) of$~u~$and $v,~$and it is denoted by$~u\ominus v$.
Furthermore, the generalized Hukuhara difference ($gH$-difference) of
two fuzzy numbers $u, v\in\mathbb{R}_{\mathcal{F}}$ is defined as
follows
\[u\ominus_{gH}v=w\Leftrightarrow \left\{\begin{array}{ll}
(i)~~u=v\oplus w,
  \\
  or
  \\
(ii)~~v=u\oplus (-1)w.
\end{array}
\right.\] It is easy to show that conditions (i) and (ii) are valid
if and only if$~w~$is a crisp number. The conditions of the
existence of $u\ominus_{gH}v\in\mathbb{R}_{\mathcal{F}}$ are given
in \cite{bs}. Through the whole of the paper, we suppose that the $gH$-difference exists.
\end{definition}
In this paper, the meaning of fuzzy-valued function is a function
$f:A\rightarrow \mathbb{R}_{\mathcal{F}},~A\in\mathbb{R}$ where
$\mathbb{R}$ is the set of all real numbers and
$[f(t)]_r=[f^{-}(t;r),f^{+}(t;r)]$ so called the $r$-cut or
parametric form of the fuzzy-valued function $f$.
\begin{definition} A fuzzy-valued function $f:[a, b]\rightarrow\mathbb{R}_{\mathcal{F}}$ is
said to be continuous at $t_0\in[a, b]$ if for each $\varepsilon>0$
there is $\delta>0$ such that $\mathcal{H}(f
(t),f(t_0))<\varepsilon$, whenever $t\in[a,b]$ and $\mid
t-t_0\mid<\delta$. We say that $f$ is fuzzy continuous on $[a,b]$ if
$f$ is continuous at each $t_0\in[a,b]$.
\end{definition}
Throughout the rest of this paper, the notation
$\mathcal{C}_f([a,b],\mathbb{R}_{\mathcal{F}})$ is called the set of
fuzzy-valued continuous functions which are defined on $[a,b]$.

If $f:[a, b]\subseteq\mathbb{R}\rightarrow\mathbb{R}_{\mathcal{F}}$
be continuous by the metric $\mathcal{H}$ then $\int^t_a f (s)ds$,
is a continuous function in $t\in[a, b]$ and the function $f$ is
integrable on $[a,b]$. Furthermore it holds
\[\bigg[\int_a^{b}f(s)ds\bigg]_r=\left[\int_a^{b}f^{-}(s;r)ds,\int_a^{b}f^{+}(s;r)ds\right].\]
\begin{definition} Let $f:[a,b]\rightarrow \mathbb{R}_{\mathcal{F}}$, $t_0\in(a,b)$ with $f^{-}(t;r)$ and $f^{+}(t;r)$ both differentiable
at $t_0$ for all $r\in[0,1]$ and $Df_{gH}$ (gH-derivative) exists:
\begin{itemize}
 \item [I.]The function $f$ is $^{F}[(i)-gH]$-differentiable at $t_0$ if
$[Df_{i.gH}(t_0)]_r=[Df^{-}(t_0;r),Df^{+}(t_0;r)]$.
\item [II.]The function $f$ is $^{F}[(ii)-gH]$-differentiable at $t_0$ if
$[Df _{ii.gH}(t_0)]_r=[Df^{+}(t_0;r),Df^{-}(t_0;r)].$
\end{itemize}
\end{definition}
\section{Definitions and Properties of Fractional $gH$-Differentiability}
In this section, let us focus on some definitions and properties
related to the fuzzy fractional generalized Hukuhara derivative
which are useful in the sequel of this paper.
\begin{definition}\cite{book}\label{d:3.1}
Let $f(t)$ be a fuzzy Lebesque integrable function. The fuzzy
Riemann-Liouville fractional (for short $(F.RL)$-fractional)
integral of order $\alpha>0$ is defined as follows
\[^{F.RL}I^{\alpha}_{[a,t]}f(t)=\frac{1}{\Gamma(\alpha)}\int^{t}_{a}(t-s)^{\alpha-1}f(s)ds.\]
\end{definition}
\begin{definition}\cite{book}\label{d:3.2}
Let $~f:[a,b]\rightarrow\mathbb{R}_{\mathcal{F}}~$. The fuzzy
fractional derivative of $f(t)$ in the Caputo sense is in the
following form
\[^{FC}D^{\alpha}_{\ast}f(t)=~^{F.RL}I^{m-\alpha}_{[a,t]}(D^{m}f_{gH})(t)=
\frac{1}{\Gamma(m-\alpha)}\int^{t}_{a}(t-s)^{m-\alpha-1}D^{m}f_{gH}(s)ds,\]
\[~m-1<\alpha<m,~m\in \mathbb{N},~t>a,\]
where $\forall m\in\mathbb{N},~D^{m}f_{gH}(s)$ ($gH$-derivatives of
$f$) are integrable. In this paper, we consider fuzzy Caputo
generalized Hukuhara derivative (for short $^{FC}[gH]$-derivative)
of order$~0<\alpha\leq1$, for fuzzy-valued function $f$, so the
$^{FC}[gH]$-derivative will be expressed by
\begin{equation}\label{eq:3.2}
^{FC}D^{\alpha}_{\ast}f(t)=~^{F.RL}I^{1-\alpha}_{[a,t]}(Df_{gH})(t)=
\frac{1}{\Gamma (1-\alpha)}\int_{a}^{t}(t-s) ^{-\alpha}Df_{gH}(s)
ds,~~~t>a.
\end{equation}
\end{definition}
\begin{lemma}\label{l:3.1}
Let $f:[a, b]\subseteq\mathbb{R}\rightarrow\mathbb{R}_{\mathcal{F}}$
be continuous. Then $^{F.RL}I_{[a,t]}^\alpha f (t)$, for
$0<\alpha\leq1$ and $t\in[a, b]$ is a continuous function.
\end{lemma}
$\mathbf{Proof}.$ Under assumptions of the continuous functions,
$f(s)$ is a fuzzy Lebesque integrable function. On the other hand,
since $\forall~0<\alpha\leq1,~(t-s)^{\alpha-1}\geq0$ is continuous,
so $\int^t_a (t-s)^{\alpha-1}f (s)ds$ is a continuous function and
as a result $^{F.RL}I_{[a,t]}^\alpha f (t)$ is a continuous function
in $t\in[a,b]$.
\begin{lemma}\label{l:3.2}
Let
$f\in\mathcal{C}_f(\mathbb{R},\mathbb{R}_{\mathcal{F}}),~m\in\mathbb{N}$.
Then the fuzzy Riemann-Liouvil fractional integrals
$^{F.RL}I^\alpha_{[a,t_{m-1}]}f(t_{m-1}),~^{F.RL}I^\alpha_{[a,t_{m-2}]}(^{F.RL}I^\alpha_{[a,t_{m-1}]}f)(t_{m-1})),...,^{F.RL}I^\alpha_{[a,t]}(^{F.RL}I^\alpha_{[a,t_{1}]}...(^{F.RL}I^\alpha_{[a,t_{m-2}]}$
\\$(^{F.RL}I^\alpha_{[a,t_{m-1}]}f)(t_{m-1}))...)$
for  $0<\alpha\leq1$, are continuous functions in
$t_{m-1},t_{m-2},...,t$, respectively. Here
$t_{m-1},t_{m-2},...,t\geq a$ and they are real numbers.
\end{lemma}
$\mathbf{Proof}$ This lemma is a fairly straightforward
generalization of Lemma $\ref{l:3.1}$. The proof will be done by
introducing on $m\in\mathbb{N}$. Assume that the lemma holds for
$(m)$-times applying operator $(F.RL)$-fractional integrating for
function $f$, we will prove it will correct for $(m+1)$-times
applying operator $(F.RL)$-fractional integrating for function $f$.
By Lemma $\ref{l:3.1}$, since
$f\in\mathcal{C}_f(\mathbb{R},\mathbb{R}_{\mathcal{F}})$ thus
$^{F.RL}I^\alpha_{[a,t_{m-1}]}f(t_{m-1})$ is a continuous function
in $t_{m-1}$. Furthermore, under the hypothesis of induction,
$$^{F.RL}I^\alpha_{[a,t_{m-2}]}(^{F.RL}I^\alpha_{[a,t_{m-1}]}f)(t_{m-1}),...,
\overbrace{^{F.RL}I^\alpha_{[a,t]}(^{F.RL}I^\alpha_{[a,t_{1}]}...(^{F.RL}I^\alpha_{[a,t_{m-2}]}
(^{F.RL}I^\alpha_{[a,t_{m-1}]}}^{(m)-times}f)(t_{m-1}))...),$$ are
continuous functions in $t_{m-1},t_{m-2},t_{m-3},...,t,$
respectively. It follows easily that
$$\overbrace{^{F.RL}I^\alpha_{[a,t_{m+1}]}(^{F.RL}I^\alpha_{[a,t]}(^{F.RL}I^\alpha_{[a,t_{1}]}...(^{F.RL}I^\alpha_{[a,t_{m-2}]}
(^{F.RL}I^\alpha_{[a,t_{m-1}]}}^{(m+1)-times}f)(t_{m-1}))...)),$$ is
a continuous function in $t_{m+1}$, which is our claim.
\begin{definition}\cite{book}\label{d:3.3}
Let$~f:[a,b]\rightarrow \mathbb{R}_{\mathcal{F}}~$be the fuzzy
Caputo generalized Hukuhara differentiable (for
short$~^{FC}[gH]$-differentiable) at$~t_{0}\in [a,b].$ Thus $~f~$is
$^{FC}[(i)-gH]$-differentiable at$~t_{0}\in [a,b]~$if for $0\leq
r\leq1$
\[~~~~[^{FC}D^{\alpha}_{\ast}f_{i.gH}(t_{0})]_r=[^{C}D^{\alpha}_{\ast}f^{-}(t_{0};r),~^{C}D^{\alpha}_{\ast}
f^{+}(t_{0};r)],\] and that $f$ is $^{FC}[(ii)-gH]$-differentiable
at $t_{0}~$if
\[~~~~[^{FC}D^{\alpha}_{\ast}f_{ii.gH}(t_{0})]_r=[^{C}D^{\alpha}_{\ast}f^{+}(t_{0};r),~^{C}D^{\alpha}_{\ast}
f^{-}(t_{0};r)],\]
where~\[^{C}D^{\alpha}_{\ast}f^{-}(t_{0};r)=\frac{1}{\Gamma(1-\alpha)}\int_{a}^{t_{0}}(t_{0}-s)
^{-\alpha}Df^{-}(s;r)ds,\]\[^{C}D^{\alpha}_{\ast}f^{+}
(t_{0};r)=\frac{1}{\Gamma(1-\alpha)}\int_{a}^{t_{0}}(t_{0}-s)
^{-\alpha}Df^{+}(s;r)ds.\]
\end{definition}
\begin{definition}\cite{book}\label{d:3.4}
Let $f:[a,b]\rightarrow\mathbb{R}_{\mathcal{F}}$ be a fuzzy-valued
function on$~[a,b]$. A point$~t_{0}\in [a,b]~$ is said to be a
switching point for the $^{FC}[gH]$-differentiability of$~f,$ if in
any neighborhood  $V~$of$~t_{0}~$there exist points
$~t_{1}<t_{0}<t_{2}~$such that
\begin{itemize}
\item [(type I)] $f$ is $^{FC}[(i)-gH]$-differentiable at $t_{1}$ while $f$ is not $^{FC}[(ii)-gH]$-differentiable at $t_{1}$, and

$ f$ is $^{FC}[(ii)-gH]$-differentiable at $t_{2}$ while $f$ is not
$^{FC}[(i)-gH]$-differentiable at $t_{2},$

or
\item [(type II)] $f$ is $^{FC}[(ii)-gH]$-differentiable at $t_{1}$ while $f$ is not $^{FC}[(i)-gH]$-differentiable at $t_{1},$

 and $f$ is $^{FC}[(i)-gH]$-differentiable at $t_{2}$ while $f$ is not $^{FC}[(ii)-gH]$-differentiable at $t_{2}.$
\end{itemize}
\end{definition}
\begin{theor}\cite{al}\label{t:3.4}
If $f:[a,b]\rightarrow \mathbb{R}_{\mathcal{F}}$,
$[f(t)]_r=[f^{-}(t;r),f^{+}(t;r)]$ and $f$ is integrable for $~0\leq
r\leq1,~t\in [a,b]$ and $\alpha, \beta>0$ then we have
\[^{F.RL}I^{\alpha}_{[a,t]}(^{F.RL}I^{\beta}_{[a,t]}f)(t)={}^{F.RL}I^{\alpha+\beta}_{[a,t]}f(t).\]
\end{theor}
\begin{lemma}\cite{al}\label{l:3.3}
Suppose that $f:[a,b]\rightarrow \mathbb{R}_{\mathcal{F}}$ be a
fuzzy-valued function and $D f_{gH}$ is exist, then for
$0<\alpha\leq1,$
\[^{F.RL}I^{\alpha}_{[a,t]}(^{FC}D^{\alpha}_{\ast}f)(t)=f(t)\ominus_{gH}f(a),~0\leq r\leq1.\]
\end{lemma}
The principal significance of this lemma is in the following
theorem:
\begin{theor}\cite{al}\label{t:3.5}
Let $f:[a,b] \rightarrow\mathbb{R}_{\mathcal{F}}$ be the  fractional
$gH$-differentiable such that type of Caputo differentiability $f$
in $[a,b]$ does not change. Then for $a\leq t\leq b$ and
$0<\alpha\leq1,$
\begin{itemize}
\item[\textbf{I.}]  If $f(s)$ is $^{FC}[(i)-gH]$-differentiable then $^{FC}D^{\alpha}_{\ast}f_{i.gH}(t)$ is $(F.RL)$-integrable over $[a,b]$ and
\[f(t)=f(a)\oplus~^{F.RL}I^{\alpha}_{[a,t]}(^{FC}D^{\alpha}_{\ast}f_{i.gH})(t),\]
\item[\textbf{II.}]  If $f(s)$ is $^{FC}[(ii)-gH]$-differentiable then $^{FC}D^{\alpha}_{\ast}f _{ii.gH}(t)$ is $(F.RL)$-integrable over $[a,b]$ and
\[f(t)=f(a)\ominus(-1)~^{F.RL}I^{\alpha}_{[a,t]}(^{FC}D^{\alpha}_{\ast}f _{ii.gH})(t).\]
\end{itemize}
\end{theor}
\begin{lemma}\label{l:3.4}
Suppose that $f:[a,b] \rightarrow\mathbb{R}_{\mathcal{F}}$ is the
fractional $gH$-differentiable and
$^{FC}D^{\alpha}_{\ast}f_{gH}(t)\in
\mathcal{C}_f([a,b],\Rb_{\mathcal{F}})$ then for $0<\alpha\leq1,$
\[^{F.RL}I_{[t,a]}^{\alpha}(^{FC}D^{\alpha}_{\ast}f_{i.gH})(t)=(-1)\odot ~^{F.RL}I_{[a,t]}^{\alpha}(^{FC}D^{\alpha}_{\ast}f _{ii.gH})(t),\]
\end{lemma}
$\mathbf{Proof.}$ Since $^{FC}D^{\alpha}_{\ast}f_{i.gH}(t)$ is
continuous, it follows that $^{FC}D^{\alpha}_{\ast}f_{i.gH}(t)$ is
the Riemann-Liouville integrable, and by using Lemma $\ref{l:3.3}$
for $0\leq r\leq 1$
\begin{eqnarray}
\nonumber [^{F.RL}I_{[t,a]}^{\alpha}(^{FC}D^{\alpha}_{\ast}f_{i.gH})(t)]_r&=&[^{F.RL}I_{[t,a]}^{\alpha}(^{FC}D^{\alpha}_{\ast}f^{-})(t;r),^{F.RL}I_{[t,a]}^{\alpha}(^{FC}D^{\alpha}_{\ast}f^{+})(t;r)]\\
\nonumber &=&[f^{-}(a;r)-f^{-}(t;r) , f^{+}(a;r)-f^{+}(t;r)]\\
\label{eq:3.3}&=&[f(a)\ominus f(t)]_r.
\end{eqnarray}
Moreover,
\begin{eqnarray}
\nonumber [^{F.RL}I_{[a,t]}^{\alpha}(^{FC}D^{\alpha}_{\ast}f _{ii.gH})(t)]_r&=&[^{F.RL}I_{[a,t]}^{\alpha}(^{FC}D^{\alpha}_{\ast}f^{+})(t;r),^{F.RL}I_{[a,t]}^{\alpha}(^{FC}D^{\alpha}_{\ast}f^{-})(t;r)]\\
\nonumber &=&[f^{+}(t;r)-f^{+}(a;r) , f^{-}(t;r)-f^{-}(a;r)]\\
 \label{eq:3.4}&=&[(-1)\odot(f(a)\ominus f(t))]_r.
\end{eqnarray}
By combining Eqs ($\ref{eq:3.3}$) with ($\ref{eq:3.4}$) the lemma is
proved.
\begin{theor}\label{t:3.6}
Let $^{FC}D^{\alpha}_{\ast}f:[a,b]
\rightarrow\mathbb{R}_{\mathcal{F}}$ and
$^{FC}D^{n\alpha}_{\ast}f\in\mathcal{C}_{f}([a,t],
\mathbb{R}_{\mathcal{F}}).$  For all $t\in [a,b]$ and
$0<\alpha\leq1,$
\begin{itemize}
\item[\textbf{I.}]Let $^{FC}D^{i\alpha}_{\ast}f,~i= 1, ..., n$ be the $^{FC}[(i)-gH]$-differentiable
and they do not change in the type of differentiability on $[a,
b],$ then
\[^{FC}D^{(i-1)\alpha}_{\ast}f_{i.gH}(t)=~^{FC}D^{(i-1)\alpha}_{\ast}f_{i.gH}(a)\oplus~^{F.RL}I^{\alpha}_{[a,t]}(^{FC}D^{i\alpha}_{\ast}f_{i.gH})(t).\]
\item[\textbf{II.}]If $^{FC}D^{i\alpha}_{\ast}f,~i= 1, ..., n$ are $^{FC}[(ii)-gH]$-differentiable and the type of
 their differentiability does not change in the interval $[a, b],$ then
\[^{FC}D^{(i-1)\alpha}_{\ast}f_{ii.gH}(t)=~^{FC}D^{(i-1)\alpha}_{\ast}f_{ii.gH}(a)\oplus~^{F.RL}I^{\alpha}_{[a,t]}(^{FC}D^{i\alpha}_{\ast}f_{ii.gH})(t).\]
\item[\textbf{III.}]Assume that  $^{FC}D^{i\alpha}_{\ast}f, i=2k-1,~k\in\mathbb{N}$ are the $^{FC}[(i)-gH]$-differentiable
and they are $^{FC}[(ii)-gH]$-differentiable, for
$i=2k,~k\in\mathbb{N}$ then
\[^{FC}D^{(i-1)\alpha}_{\ast}f_{i.gH}(t)=~^{FC}D^{(i-1)\alpha}_{\ast}f_{i.gH}(a)\ominus(-1)^{F.RL}I^{\alpha}_{[a,t]}(^{FC}D^{i\alpha}_{\ast}f_{ii.gH})(t).\]
\item[\textbf{IV.}]Suppose that $^{FC}D^{i\alpha}_{\ast}f,~i=2k-1,~k\in\mathbb{N}$ are $^{FC}[(ii)-gH]$-differentiable
and they are $^{FC}[(i)-gH]$-differentiable for
$i=2k,~k\in\mathbb{N}$, so
\[^{FC}D^{(i-1)\alpha}_{\ast}f_{ii.gH}(t)=~^{FC}D^{(i-1)\alpha}_{\ast}f_{ii.gH}(a)\ominus(-1)^{F.RL}I^{\alpha}_{[a,t]}(^{FC}D^{i\alpha}_{\ast}f_{i.gH})(t).\]
\end{itemize}
\end{theor}
$\mathbf{Proof.}$ By assuming $^{FC}D^{i\alpha}_{\ast}f \in
\mathcal{C}_{f}([a,b],\mathbb{R}_{\mathcal{F}}),~i=0,...,n$ we give
the proof only for parts \textbf{II} and \textbf{III}. Proving the
other parts are similar.

$\textbf{II.}$ Our proof starts with the observation that
$^{FC}D^{i\alpha}_{\ast}f,~i=1,...,n$ are
$^{FC}[(ii)-gH]$-differentiable. Hence, using properties of fuzzy
Caputo derivative and Theorem $\ref{t:3.5}$, we have
\begin{align*}
&[^{F.RL}I^{\alpha}_{[a,t]}(^{FC}D^{i\alpha}_{\ast}f_{ii.gH})(t)]_r\\
&=[^{RL}I^{\alpha}_{[a,t]}(^{C}D^{i\alpha}_{\ast}f^{+})(t;r),~^{RL}I^{\alpha}_{[a,t]}(^{C}D^{i\alpha}_{\ast}f^{-})(t;r)]\\
&=[^{RL}I^{\alpha}_{[a,t]}.^{C}D^{\alpha}_{\ast}(^{C}D^{(i-1)\alpha}_{\ast}f^{+})(t;r),~^{RL}I^{\alpha}_{[a,t]}.^{C}D^{\alpha}_{[a,t]}(^{C}D^{(i-1)\alpha}_{\ast}f^{-})(t;r)]\\
&=[^{C}D^{(i-1)\alpha}_{\ast}f^{+}(t;r)-~^{C}D^{(i-1)\alpha}_{\ast}f^{+}(a;r),~^{C}D^{(i-1)\alpha}_{\ast}f^{-}(t;r)-~^{C}D^{(i-1)\alpha}_{\ast}f^{-}(a;r)]\\
&=[^{C}D^{(i-1)\alpha}_{\ast}f^{+}(t;r),~^{C}D^{(i-1)\alpha}_{\ast}f^{-}(t;r)]-[^{C}D^{(i-1)\alpha}_{\ast}f^{+}(a;r),~^{C}D^{(i-1)\alpha}_{\ast}f^{-}(a;r)]\\
&=[^{FC}D^{(i-1)\alpha}_{\ast}f_{ii.gH}(t)\ominus{}^{FC}D^{(i-1)\alpha}_{\ast}f_{ii.gH}(a)]_r.
\end{align*}
Thus, we obtain
\[^{FC}D^{(i-1)\alpha}_{\ast}f_{ii.gH}(t)=~^{FC}D^{(i-1)\alpha}_{\ast}f_{ii.gH}(a)\oplus ^{F.RL}I^{\alpha}_{[a,t]}(^{FC}D^{i\alpha}_{\ast}f_{ii.gH})(t).\]

$\textbf{III.}$ Under the conditions stated in the part III,
$^{FC}D^{i\alpha}_{\ast}f$ is $^{FC}[(i)-gH]$-differentiable for
$i=2k-1,~k\in\mathbb{N}$ and it is $^{FC}[(ii)-gH]$-differentiable
for $i=2k,~k\in\mathbb{N}$. In the sense of Section 2 and by Theorem
$\ref{t:3.5}$, we get
\begin{align*}
&[^{FC}D^{(i-1)\alpha}_{\ast}f_{i.gH}(t)\oplus(-1)^{F.RL}I^{\alpha}_{[a,t]}(^{FC}D^{i\alpha}_{\ast}f_{ii.gH})(t)]_r\\
&=[^{C}D^{(i-1)\alpha}_{\ast}f^{-}(t;r),^{C}D^{(i-1)\alpha}_{\ast}f_{+}(t;r)]+[-^{RL}I^{\alpha}_{[a,t]}(^{C}D^{i\alpha}_{\ast}f^{-})(t;r),-^{RL}I^{\alpha}_{[a,t]}(^{C}D^{i\alpha}_{\ast}f^{+})(t;r)]\\
&=[^{C}D^{(i-1)\alpha}_{\ast}f^{-}(t;r),~^{C}D^{(i-1)\alpha}_{\ast}f^{+}(t;r)]\\
&+[^{C}D^{(i-1)\alpha}_{\ast}f^{-}(a;r)-~^{C}D^{(i-1)\alpha}_{\ast}f^{-}(t;r),~^{C}D^{(i-1)\alpha}_{\ast}f^{+}(a;r)-~^{C}D^{(i-1)\alpha}_{\ast}f^{+}(t;r)]\\
&=[^{C}D^{(i-1)\alpha}_{\ast}f^{-}(a;r),~^{C}D^{(i-1)\alpha}_{\ast}f^{+}(a;r)]=[^{FC}D^{(i-1)\alpha}_{\ast}f_{i.gH}(a)]_r,
\end{align*}
which completes the proof.
\section{Fuzzy Generalized Taylor Theorem}
\begin{theor}\label{t:3.7}
Let $T=[a,a+\beta]\subset \mathbb{R}$, with $\beta>0$ and $^{FC}D^{i\alpha}_{\ast}f \in\mathcal{C}_f([a,b],\mathbb{R}_{f}),~i=1,...,n$. For $t\in T,~ 0<\alpha\leq 1$\\

$\textbf{I.}$ If $~^{FC}D^{i\alpha}_{\ast}f,~i=0,1,...,n-1$ are
$^{FC}[(i)-gH]$-differentiable, provided that type of fuzzy Caputo
differentiability has no change. Then
\begin{eqnarray}
\nonumber f(t)&=&f(a)\oplus ~^{FC}D^{\alpha}_{\ast}f_{i.gH} (a)\odot\frac{(t-a)^{\alpha}}{\Gamma(\alpha+1)}\oplus~^{FC}D^{2\alpha}_{\ast}f_{i.gH} (a)\odot \frac{(t-a)^{2\alpha}}{\Gamma (2\alpha+1)}\\
\nonumber &\oplus &...\oplus~^{FC}D^{(n-1)\alpha}_{\ast}f_{i.gH}
(a)\odot \frac{(t-a)^{(n-1)\alpha}}{\Gamma ((n-1)\alpha+1)}\oplus
R_{n}(a,t),
\end{eqnarray}
where
$R_{n}(a,t):=~^{F.RL}I^{\alpha}_{[a,t]}(^{F.RL}I^{\alpha}_{[a,t_{1}]}
...(^{F.RL}I^{\alpha}_{[a,t_{n-1}]}(^{FC}D^{n\alpha}_{\ast}f_{i.gH})(t_{n}))...).$\\

$\textbf{II.}$ If $~^{FC}D^{i\alpha}_{\ast}f,~i=0,1,...,n-1$ are
$^{FC}[(ii)-gH]$-differentiable, provided that type of fuzzy Caputo
differentiability has no change. Then
\begin{eqnarray}
\nonumber f(t)&=&f(a)\ominus(-1) ~^{FC}D^{\alpha}_{\ast}f _{ii.gH}(a)\odot\frac{(t-a)^{\alpha}}{\Gamma(\alpha+1)}\ominus(-1)~^{FC}D^{2\alpha}_{\ast}f _{ii.gH} (a)\odot\frac{(t-a)^{2\alpha}}{\Gamma (2\alpha+1)}\\
\nonumber &\ominus &(-1) ...\ominus(-1)~^{FC}D^{(n-1)\alpha}_{\ast}f
_{ii.gH} (a)\odot \frac{(t-a)^{(n-1)\alpha}}{\Gamma
((n-1)\alpha+1)}\ominus(-1)R   _{n}(a,t),
\end{eqnarray}
where
$R_{n}(a,t):=~^{F.RL}I^{\alpha}_{[a,t]}(^{F.RL}I^{\alpha}_{[a,t_{1}]}
...(^{F.RL}I^{\alpha}_{[a,t_{n-1}]}(^{FC}D^{n\alpha}_{\ast}f _{ii.gH})(t_{n}))...).$\\

$\textbf{III.}$ If $~^{FC}D^{i\alpha}_{\ast}f,~i=2k-1,~k\in
\mathbb{N}$ are $^{FC}[(i)-gH]$-differentiable and
$^{FC}D^{i\alpha}_{\ast}f,~i=2k,~k\in \mathbb{N}\cup \{0\}$ are
$^{FC}[(ii)-gH]$-differentiable, then
\begin{eqnarray}
\nonumber f(t)&=&f(a)\ominus(-1)~^{FC}D^{\alpha}_{\ast}f _{ii.gH} (a)\odot\frac{(t-a)^{\alpha}}{\Gamma(\alpha+1)}\oplus~^{FC}D^{2\alpha}_{\ast}f_{i.gH} (a)\odot\frac{(t-a)^{2\alpha}}{\Gamma (2\alpha+1)}\\
\nonumber &\ominus &(-1) ...\ominus(-1)~^{FC}D^{(\frac{i}{2}-1)\alpha}_{\ast}f _{ii.gH} (a)\odot \frac{(t-a)^{(\frac{i}{2}-1)\alpha}}{\Gamma (\frac{i}{2}\alpha)}\\
\nonumber &\oplus &~^{FC}D^{(\frac{i}{2})\alpha}_{\ast}f_{i.gH}
(a)\odot \frac{(t-a)^{(\frac{i}{2})\alpha}}{\Gamma
(\frac{i}{2}\alpha+1)}\ominus(-1)...\ominus(-1)R_{n}(a,t),
\end{eqnarray}
where
$R_{n}(a,t):=~^{F.RL}I^{\alpha}_{[a,t]}(^{F.RL}I^{\alpha}_{[a,t_{1}]}
...(^{F.RL}I^{\alpha}_{[a,t_{n-1}]}(^{FC}D^{n\alpha}_{\ast}f_{i.gH})(t_{n}))...).$\\

$\textbf{IV.}$ For $~^{FC}D^{n\alpha}_{\ast}f
\in\mathcal{C}_{f}([a,b],\mathbb{R}_{f}),~n\geq3$, suppose that $f$
on $[a, \xi]$ is $^{FC}[(ii) -[gH]$-differentiable and on $[\xi, b]$
is $^{FC}[(i)-gH]$-differentiable, in fact $\xi$ is switching point
(type II) for $\alpha$-order derivative of $f$. Moreover, for
$t_0\in[a, \xi]$, let $2\alpha$-order derivative of $f$ in $\xi_{1}$
of $[t_{0}, \xi]$ have switching point (type I). On the other hand,
the type of differentiability for $^{FC}D^{i\alpha}_{\ast}f, i\leq
n$ on $[\xi, b]$ does not change. So
\begin{eqnarray}
\nonumber f(t)&=&f(t_0)\ominus(-1)^{FC}D^{\alpha}_{\ast}f_{ii.gH}(t_0)\odot \frac{(\xi-t_0)^{\alpha}}{\Gamma(\alpha+1)}\oplus (-1)^{FC}D^{2\alpha}_{\ast}f_{i.gH}(t_0)\odot \frac{(t_0-\xi_{1})^{\alpha}}{\Gamma(\alpha+1)}\\
\nonumber &\odot &\frac{(\xi -t_0)^{\alpha}}{\Gamma(\alpha+1)}\ominus(-1)^{FC}D^{2\alpha}_{\ast}f_{ii.gH}(\xi_{1})\odot\bigg[\frac{(\xi -\xi_{1})^{2\alpha}}{\Gamma(2\alpha+1)}-\frac{(t_0-\xi_{1})^{2\alpha}}{\Gamma(2\alpha+1)}\bigg]\\
\nonumber &\oplus &^{FC}D^{\alpha}_{\ast}f_{i.gH
}(\xi)\odot\frac{(t-\xi)^{\alpha}}{\Gamma(\alpha+1)}\oplus ~^{FC}D^{2\alpha}_{\ast}f_{i.gH}(\xi)\odot \frac{(t-\xi)^{2\alpha}}{\Gamma(2\alpha+1)}\\
\nonumber &\oplus &^{F.RL}I_{[t_0,\xi]}^{\alpha}.^{F.RL}I_{[t_0,\xi_{1}]}.^{F.RL}I_{[t_0,t_2]}^{\alpha}(^{FC}D^{3\alpha}_{\ast}f_{i.gH})(t_4)\\
\nonumber &\ominus &(-1)^{F.RL}I_{[t_0,\xi]}^{\alpha}.^{F.RL}I_{[\xi_{1},t_1]}.^{F.RL}I_{[\xi_{1},t_3]}^{\alpha}(^{FC}D^{3\alpha}_{\ast}f_{ii.gH})(t_5).\\
\nonumber &\oplus &~^{F.RL}I_{[\xi,t]}.^{F.RL}I_{[\xi
,s_1]}^{\alpha}.^{F.RL}I_{[a,s_2]}^{\alpha}(^{FC}D^{3\alpha}_{\ast}f_{i.gH})(s_3).
\end{eqnarray}
$\mathbf{Proof.}$ Under the assumptions that
$^{FC}D^{i\alpha}_{\ast}f
\in\mathcal{C}_f([a,b],\mathbb{R}_{f}),~i=1,...,n,$ we conclude that
$^{FC}D^{i\alpha}_{\ast}f$ are $(F.RL)$-fractional integrable on
$T$,

$\textbf{I.}$ Since $f$ is a continuous function and
$^{FC}[(i)-gH]$-differentiable, by Theorem $\ref{t:3.5}$, we get
\[f(t)=f(a)\oplus~^{F.RL}I^{\alpha}_{[a,t]}(^{FC}D^{\alpha}_{\ast}f_{i.gH})(t_1),\]
and Theorem $\ref{t:3.6}$ yields
\[{}^{FC}D^{\alpha}_{\ast}f_{i.gH}(t_1)={}^{FC}D^{\alpha}_{\ast}f_{i.gH}(a)\oplus ~^{F.RL}I^{\alpha}_{[a,t_{1}]}(^{FC}D^{2\alpha}_{\ast}f_{i.gH})(t_{2}).\]
Therefore
\begin{align*}
^{F.RL}I^{\alpha}_{[a,t]}(^{FC}D^{\alpha}_{\ast}f_{i.gH})(t_1)={}^{FC}D^{\alpha}_{\ast}f_{i.gH}(a)\odot\frac{(t-a)^{\alpha}}{\Gamma(\alpha+1)}\oplus
~^{F.RL}I^{\alpha}_{[a,t]}.~^{F.RL}I^{\alpha}_{[a,t_{1}]}(^{FC}D^{2\alpha}_{\ast}f_{i.gH})(t_{2}).
\end{align*}
Since the last double $(F.RL)$-fractional integral belongs to
$\mathbb{R}_f$ and by using Lemma $\ref{l:3.2}$ we have
\[f(t)=f(a)\oplus{}^{FC}D^{\alpha}_{\ast}f_{i.gH}(a)\odot\frac{(t-a)^{\alpha}}{\Gamma(\alpha+1)}\oplus ~^{F.RL}I^{\alpha}_{[a,t]}.~^{F.RL}I^{\alpha}_{[a,t_{1}]}(^{FC}D^{2\alpha}_{\ast}f_{i.gH})(t_{2}).\]
By similar argument,
\[{}^{FC}D^{2\alpha}_{\ast}f_{i.gH}(t_{2})={}^{FC}D^{2\alpha}_{\ast}f_{i.gH}(a)\oplus{}^{F.RL}I^{\alpha}_{[a,t_{2}]}(^{FC}D^{3\alpha}_{\ast}f_{i.gH})(t_{3}).\]
Applying operator $(F.RL)$-fractional integral to
$(^{FC}D^{2\alpha}_{\ast}f_{i.gH})(t_{2})$, we obtain
\begin{align*}
^{F.RL}I^{\alpha}_{[a,t_{1}]}(^{FC}D^{2\alpha}_{\ast}f_{i.gH})(t_{2})
=^{F.RL}I^{\alpha}_{[a,t_{1}]}(^{FC}D^{2\alpha}_{\ast}f_{i.gH})(a)\oplus
~^{F.RL}I^{\alpha}_{[a,t_{1}]}.~^{F.RL}I^{\alpha}_{[a,t_{2}]}(^{FC}D^{3\alpha}_{\ast}f_{i.gH})(t_{3}),
\end{align*}
thus
\begin{align*}
^{F.RL}I^{\alpha}_{[a,t_{1}]}(^{FC}D^{2\alpha}_{\ast}f_{i.gH})(t_{2})={}^{FC}D^{2\alpha}_{\ast}f_{i.gH}(a)\odot\frac{(t_{1}-a)^{\alpha}}{\Gamma(\alpha+1)}\oplus
^{F.RL}I^{\alpha}_{[a,t_{1}]}.~^{F.RL}I^{\alpha}_{[a,t_{2}]}(^{FC}D^{3\alpha}_{\ast}f_{i.gH})(t_{3}),
\end{align*}
furthermore
\begin{eqnarray}
\nonumber ^{F.RL}I^{\alpha}_{[a,t]}.^{F.RL}I^{\alpha}_{[a,t_{1}]}(^{FC}D^{2\alpha}_{\ast}f_{i.gH})(t_{2})&=&{}^{FC}D^{2\alpha}_{\ast}f_{i.gH}(a)\odot \frac{(t-a)^{2\alpha}}{\Gamma(2\alpha+1)}\\
\nonumber &\oplus
&^{F.RL}I^{\alpha}_{[a,t]}.~^{F.RL}I^{\alpha}_{[a,t_{1}]}.~^{F.RL}I^{\alpha}_{[a,t_{2}]}(^{FC}D^{3\alpha}_{\ast}f_{i.gH})(t_{3}).
\end{eqnarray}
The last triple integral belongs to $\mathbb{R}_f$. By Lemma
$\ref{l:3.2}$ we get
\begin{eqnarray}
\nonumber f(t)&=&f(a)\oplus{}^{FC}D^{\alpha}_{\ast}f_{i.gH}(a)\odot\frac{(t-a)^{\alpha}}{\Gamma(\alpha+1)}\oplus {}^{FC}D^{2\alpha}_{\ast}f_{i.gH}(a)\odot \frac{(t-a)^{2\alpha}}{\Gamma(2\alpha+1)}\\
\nonumber &\oplus
&^{F.RL}I^{\alpha}_{[a,t]}.~^{F.RL}I^{\alpha}_{[a,t_{1}]}.~^{F.RL}I^{\alpha}_{[a,t_{2}]}(^{FC}D^{3\alpha}_{\ast}f_{i.gH})(t_{3}).
\end{eqnarray}
The high order of the last formula by Lemma $\ref{l:3.2}$ is a
continuous function in terms of $t$ so it belongs to $\Rb_f$. With
the same manner, we can demonstrate that part $\textbf{I}$ is
satisfied.

$\textbf{II.}$ Let $f$ is $^{FC}[(ii)-gH]$-differentiable, we have
\[f(t)=f(a)\ominus (-1)~^{F.RL}I^{\alpha}_{[a,t]}(^{FC}D^{\alpha}_{\ast}f _{ii.gH})(t_1).\]
Under the hypotheses of Theorem, type of differentiability does not
change, so by Theorem $\ref{t:3.6}$ and by attention to
$(F.RL)$-integrability of $^{FC}D^{\alpha}_{\ast}f _{ii.gH}$ on $T$
, we obtain
\[{}^{FC}D^{\alpha}_{\ast}f _{ii.gH}(t_{1})={}^{FC}D^{\alpha}_{\ast}f _{ii.gH}(a)\oplus{}^{F.RL}I^{\alpha}_{[a,t_{1}]}(^{FC}D^{2\alpha}_{\ast}f _{ii.gH})(t_{2}).\]
Applying operator $^{F.RL}I^{\alpha}_{[a,t]}$ to
$^{FC}D^{\alpha}_{\ast}f _{ii.gH}(t_1)$, gives
\begin{align*}
^{F.RL}I^{\alpha}_{[a,t]}(^{FC}D^{\alpha}_{\ast}f _{ii.gH})(t_{1})
={}^{FC}D^{\alpha}_{\ast}f
_{ii.gH}(a)\odot\frac{(t-a)^{\alpha}}{\Gamma(\alpha+1)}\oplus
~^{F.RL}I^{\alpha}_{[a,t]}.~^{F.RL}I^{\alpha}_{[a,t_{1}]}(^{FC}D^{2\alpha}_{\ast}f
_{ii.gH})(t_{2}).
\end{align*}
Lemma $\ref{l:3.2}$ implies that the last double $(F.RL)$-fractional
integral belongs to $\mathbb{R}_{f}$. So
\begin{equation}\label{eq:4.5}
f(t)=f(a)\ominus(-1){}^{FC}D^{\alpha}_{\ast}f
_{ii.gH}(a)\odot\frac{(t-a)^{\alpha}}{\Gamma(\alpha+1)}\ominus(-1)
~^{F.RL}I^{\alpha}_{[a,t]}.~^{F.RL}I^{\alpha}_{[a,t_{1}]}(^{FC}D^{2\alpha}_{\ast}f
_{ii.gH})(t_{2}).
\end{equation}
By repeating the above argument, we get
\[{}^{FC}D^{2\alpha}_{\ast}f _{ii.gH}(t_{2})={}^{FC}D^{2\alpha}_{\ast}f _{ii.gH}(a)\oplus ~^{F.RL}I^{\alpha}_{[a,t_{2}]}(^{FC}D^{3\alpha}_{\ast}f _{ii.gH})(t_{3}).\]
Therefore, we find that
\[^{F.RL}I^{\alpha}_{[a,t_{1}]}(^{FC}D^{2\alpha}_{\ast}f _{ii.gH})(t_{2})={}^{FC}D^{2\alpha}_{\ast}f _{ii.gH}(a)\odot\frac{(t_{1}-a)^{\alpha}}{\Gamma(\alpha+1)}\\
\oplus~^{F.RL}I^{\alpha}_{[a,t_{1}]}.~^{F.RL}I^{\alpha}_{[a,t_{2}]}(^{FC}D^{3\alpha}_{\ast}f
_{ii.gH})(t_{3}).\] Moreover
\begin{eqnarray}
\nonumber ^{F.RL}I^{\alpha}_{[a,t]}.^{F.RL}I^{\alpha}_{[a,t_{1}]}(^{FC}D^{2\alpha}_{\ast}f _{ii.gH})(t_{2})&=&{}^{FC}D^{2\alpha}_{\ast}f _{ii.gH}(a)\odot \frac{(t-a)^{2\alpha}}{\Gamma(2\alpha+1)}\\
\nonumber &\oplus
&{}^{F.RL}I^{\alpha}_{[a,t]}.~^{F.RL}I^{\alpha}_{[a,t_{1}]}.~^{F.RL}I^{\alpha}_{[a,t_{2}]}(^{FC}D^{3\alpha}_{\ast}f
_{ii.gH})(t_{3}).
\end{eqnarray}
By Lemma $\ref{l:3.2}$, the last triple $(F.RL)$-fractional integral
belongs to $\mathbb{R}_f$. Therefore, substituting above equation
into Eq. $(\ref{eq:4.5})$, we find that
\begin{eqnarray}
\nonumber f(t)&=&f(a)\ominus(-1){}^{FC}D^{\alpha}_{\ast}f _{ii.gH}(a)\odot\frac{(t-a)^{\alpha}}{\Gamma(\alpha+1)}\ominus(-1) {}^{FC}D^{2\alpha}_{\ast}f _{ii.gH}(a)\\
\nonumber &\odot &\frac{(t-a)^{2\alpha}}{\Gamma(2\alpha+1)}\ominus~
(-1)^{F.RL}I^{\alpha}_{[a,t]}.~^{F.RL}I^{\alpha}_{[a,t_{1}]}.~^{F.RL}I^{\alpha}_{[a,t_{2}]}(^{FC}D^{3\alpha}_{\ast}f
_{ii.gH})(t_{3}).
\end{eqnarray}
The rest of the proof runs as before.

$\textbf{III.}$ Suppose that $f$ is $^{FC}[(ii)-gH]$-differentiable.
Using Theorem $\ref{t:3.5}$ we have
\[f(t)=f(a)\ominus (-1)~^{F.RL}I^{\alpha}_{[a,t]}(^{FC}D^{\alpha}_{\ast}f _{ii.gH})(t_1).\]
Under the hypothesis of theorem, since $f$ is
$^{FC}[(ii)-gH]$-differentiable, $^{FC}D^{\alpha}_{\ast}f$ is
$^{FC}[(i)-gH]$-differentiable. So, by Theorem $\ref{t:3.6}$ we get
\[{}^{FC}D^{\alpha}_{\ast}f _{ii.gH}(t_{1})={}^{FC}D^{\alpha}_{\ast}f _{ii.gH}(a)\ominus (-1)~^{F.RL}I^{\alpha}_{[a,t_{1}]}(^{FC}D^{2\alpha}_{\ast}f_{i.gH})(t_{2}).\]
Now, applying operator $(F.RL)$-integral to
$(^{FC}D^{\alpha}_{\ast}f _{ii.gH})(t_1)$ gives
\begin{eqnarray}
\nonumber &^{F.RL}I^{\alpha}_{[a,t]}&(^{FC}D^{\alpha}_{\ast}f _{ii.gH})(t_{1})\\
\nonumber &&=~^{F.RL}I^{\alpha}_{[a,t]}(^{FC}D^{\alpha}_{\ast}f _{ii.gH})(a)\ominus (-1) ~^{F.RL}I^{\alpha}_{[a,t]}.~^{F.RL}I^{\alpha}_{[a,t_{1}]}(^{FC}D^{2\alpha}_{\ast}f_{i.gH})(t_{2})\\
\nonumber &&={}^{FC}D^{\alpha}_{\ast}f
_{ii.gH}(a)\odot\frac{(t-a)^{\alpha}}{\Gamma(\alpha+1)}\ominus (-1)
^{F.RL}I^{\alpha}_{[a,t]}.~^{F.RL}I^{\alpha}_{[a,t_{1}]}(^{FC}D^{2\alpha}_{\ast}f_{i.gH})(t_{2}).
\end{eqnarray}
Lemma $\ref{l:3.2}$ now leads to the last double (F.RL)-fractional
integral belongs to $\mathbb{R}_{f}$. So
\[f(t)=f(a)\ominus(-1){}^{FC}D^{\alpha}_{\ast}f _{ii.gH}(a)\odot\frac{(t-a)^{\alpha}}{\Gamma(\alpha+1)}\oplus~ ~^{F.RL}I^{\alpha}_{[a,t]}.~^{F.RL}I^{\alpha}_{[a,t_{1}]}(^{FC}D^{2\alpha}_{\ast}f _{i.gH})(t_{2}).\]
Similarly, since $^{FC}D^{\alpha}_{\ast}f$ is
$^{FC}[(i)-gH]$-differentiable, $^{FC}D^{2\alpha}_{\ast}f$ is
$^{FC}[(ii)-gH]$-differentiable and we get
\[{}^{FC}D^{2\alpha}_{\ast}f_{i.gH}(t_{2})={}^{FC}D^{2\alpha}_{\ast}f_{i.gH})(a)\ominus (-1)~^{F.RL}I^{\alpha}_{[a,t_{2}]}(^{FC}D^{3\alpha}_{\ast}f _{ii.gH})(t_{3}).\]
Thus
\begin{eqnarray}
\nonumber ^{F.RL}I^{\alpha}_{[a,t_{1}]}(^{FC}D^{2\alpha}_{\ast}f_{i.gH})(t_{2})&=&{}^{FC}D^{2\alpha}_{\ast}f_{i.gH}(a)\odot\frac{(t_{1}-a)^{\alpha}}{\Gamma(\alpha+1)}\\
\nonumber &\ominus
&(-1)~^{F.RL}I^{\alpha}_{[a,t_{1}]}.~^{F.RL}I^{\alpha}_{[a,t_{2}]}(^{FC}D^{3\alpha}_{\ast}f
_{ii.gH})(t_{3}).
\end{eqnarray}
Now, applying operator $^{F.RL}I^{\alpha}_{[a,t]}$ gives
\begin{eqnarray}
\nonumber ^{F.RL}I^{\alpha}_{[a,t]}.^{F.RL}I^{\alpha}_{[a,t_{1}]}(^{FC}D^{2\alpha}_{\ast}f_{i.gH})(t_{2})&=&{}^{FC}D^{2\alpha}_{\ast}f_{i.gH}(a)\odot \frac{(t-a)^{2\alpha}}{\Gamma(2\alpha+1)}\\
\nonumber &\ominus
&(-1)~^{F.RL}I^{\alpha}_{[a,t]}.~^{F.RL}I^{\alpha}_{[a,t_{1}]}.~^{F.RL}I^{\alpha}_{[a,t_{2}]}(^{FC}D^{3\alpha}_{\ast}f_{ii.gH})(t_{3}).
\end{eqnarray}
Since satisfies all the other conditions for the Lemma
$\ref{l:3.2}$, the last triple $(F.RL)$-fractional integral belongs
to $\mathbb{R}_f$. Then
\begin{eqnarray}
\nonumber f(t)&=&f(a)\ominus(-1){}^{FC}D^{\alpha}_{\ast}f_{ii.gH}(a)\odot \frac{(t-a)^{\alpha}}{\Gamma(\alpha+1)}\oplus~ {}^{FC}D^{2\alpha}_{\ast}f_{i.gH}(a)\\
\nonumber &\odot
&\frac{(t-a)^{2\alpha}}{\Gamma(2\alpha+1)}\ominus(-1)~^{F.RL}I^{\alpha}_{[a,t]}.~^{F.RL}I^{\alpha}_{[a,t_{1}]}.~^{F.RL}I^{\alpha}_{[a,t_{2}]}(^{FC}D^{3\alpha}_{\ast}f_{ii.gH})(t_{3}).
\end{eqnarray}
with  simple and  similar method, the proof for this type of
differentiability will be completed.

$\textbf{IV.}$ Since $f$ is $^{FC}[(ii)-gH]$-differentiable in
$[t_0,\xi]$, Theorem $\ref{t:3.5}$ leads to
\begin{equation}\label{eq:4.6}
f(\xi)=f(t_0)\ominus(-1)^{F.RL}I_{[t_0,\xi]}^{\alpha}(^{FC}D^{\alpha}_{\ast}f_{ii.gH})(t_1),
\end{equation}
and in the interval $[\xi,b]$, $f$ is
$^{FC}[(i)-gH]$-differentiable, so for $t\in [\xi,b]$
\begin{equation}\label{eq:4.7}
f(t)=f(\xi)\oplus~^{F.RL}I_{[\xi,t]}^{\alpha}(^{FC}D^{\alpha}_{\ast}f_{i.gH})(s_1).
\end{equation}
According to the hypothesis, we know that $\xi$ is a switching point
for differentiability $f$, thus by substituting Eq. $(\ref{eq:4.6})$
into Eq. $(\ref{eq:4.7})$ we obtain
\begin{equation}\label{eq:4.8}
f(t)=f(t_0)\ominus(-1)^{F.RL}I_{[t_0,\xi]}^{\alpha}(^{FC}D^{\alpha}_{\ast}f_{ii.gH})(t_1)\oplus~^{F.RL}I_{[\xi,t]}^{\alpha}(^{FC}D^{\alpha}_{\ast}f_{i.gH})(s_1).
\end{equation}
Consider the first $(F.RL)$-fractional integral on the right side of the Eq. $(\ref{eq:4.8})$:\\
By noting the hypothesis of theorem, the fuzzy Caputo derivative of
the function $f$ has the switching point $\xi_1$ of type I. So,
$^{FC}D^{\alpha}_{\ast}f_{ii.gH}$ is $^{FC}[(i)-gH]$-differentiable
on $[t_0,\xi_1]$, then type of differentiability can be changed. By
these conditions, the Theorem $\ref{t:3.6}$, admits that
\begin{equation}\label{eq:4.9}
^{FC}D^{\alpha}_{\ast}f_{ii.gH}(\xi_1)=~^{FC}D^{\alpha}_{\ast}f_{ii.gH}(t_0)\ominus(-1)^{F.RL}I_{[t_0
,\xi_1]}^{\alpha}(^{FC}D^{2\alpha}_{\ast}f_{i.gH})(t_2).
\end{equation}
On the other hand, we know that $^{FC}D^{\alpha}_{\ast}f_{ii.gH}$ is
$^{FC}[(ii)-gH]$-differentiable on $[\xi_1 ,\xi]$ and the type of
differentiability does not change. Thus, for $t_1\in[\xi_1 ,\xi]$
from Theorem $\ref{t:3.6}$, it follows that
\begin{equation}\label{eq:4.10}
^{FC}D^{\alpha}_{\ast}f_{ii.gH}(t_1)=~^{FC}D^{\alpha}_{\ast}f_{ii.gH}(\xi_1)\oplus~^{F.RL}I_{[\xi_{1},t_1]}^{\alpha}(^{FC}D^{2\alpha}_{\ast}f_{ii.gH})(t_3).
\end{equation}
Substituting Eq. $(\ref{eq:4.9})$ into Eq. $(\ref{eq:4.10})$ gives
\begin{eqnarray}
\nonumber ^{FC}D^{\alpha}_{\ast}f_{ii.gH}(t_1)&=&^{FC}D^{\alpha}_{\ast}f_{ii.gH}(t_0)\ominus(-1)^{F.RL}I_{[t_0 ,\xi_1]}^{\alpha}(^{FC}D^{2\alpha}_{\ast}f_{i.gH})(t_2)\\
\label{eq:4.11}&\oplus
&^{F.RL}I_{[\xi_{1},t_1]}^{\alpha}(^{FC}D^{2\alpha}_{\ast}f_{ii.gH})(t_3),
\end{eqnarray}
that
\begin{eqnarray}
\nonumber &&^{FC}D^{2\alpha}_{\ast}f_{i.gH}(t_2)=~^{FC}D^{2\alpha}_{\ast}f_{i.gH}(t_0)\oplus~^{F.RL}I_{[t_0,t_2]}^{\alpha}(^{FC}D^{3\alpha}_{\ast}f_{i.gH})(t_4),\\
\nonumber &\Rightarrow &^{F.RL}I_{[t_0,\xi_{1}]}^{\alpha}(^{FC}D^{2\alpha}_{\ast}f_{i.gH})(t_2)\\
\label{eq:4.12}&=&^{FC}D^{2\alpha}_{\ast}f_{i.gH}(t_0)\odot
\frac{(\xi_1-t_0)^{\alpha}}{\Gamma(\alpha+1)}\oplus
~^{F.RL}I_{[t_0,\xi_{1}]}.^{F.RL}I_{[t_0,t_2]}^{\alpha}(^{FC}D^{3\alpha}_{\ast}f_{i.gH})(t_4),
\end{eqnarray}
follows from Theorem $\ref{t:3.6}$ and also
\begin{eqnarray}
\nonumber &&^{FC}D^{2\alpha}_{\ast}f_{ii.gH}(t_3)=~^{FC}D^{2\alpha}_{\ast}f_{ii.gH}(\xi_{1})\oplus ~^{F.RL}I_{[\xi_{1},t_3]}^{\alpha}(^{FC}D^{3\alpha}_{\ast}f_{ii.gH})(t_5),\\
\nonumber &\Rightarrow &^{F.RL}I_{[\xi_{1},t_1]}^{\alpha}(^{FC}D^{2\alpha}_{\ast}f_{ii.gH})(t_3)\\
\label{eq:4.13}&=&^{FC}D^{2\alpha}_{\ast}f_{ii.gH}(\xi_{1})\odot
\frac{(t_1-\xi_{1})^{\alpha}}{\Gamma(\alpha+1)}\oplus
~^{F.RL}I_{[\xi_{1},t_1]}.^{F.RL}I_{[\xi_{1},t_3]}^{\alpha}(^{FC}D^{3\alpha}_{\ast}f_{ii.gH})(t_5).
\end{eqnarray}
The insertion of the Eqs. $(\ref{eq:4.12})$ and $(\ref{eq:4.13})$,
in Eq. $(\ref{eq:4.11})$ leads to obtain
\begin{eqnarray}
\nonumber ^{FC}D^{\alpha}_{\ast}f_{ii.gH}(t_1)&=&^{FC}D^{\alpha}_{\ast}f_{ii.gH}(t_0)\ominus ~^{FC}D^{2\alpha}_{\ast}f_{i.gH}(t_0)\odot \frac{(t_0-\xi_{1})^{\alpha}}{\Gamma(\alpha+1)}\oplus ~^{FC}D^{2\alpha}_{\ast}f_{ii.gH}(\xi_{1})\\
\nonumber &\odot &\frac{(t_1-\xi_{1})^{\alpha}}{\Gamma(\alpha+1)}\ominus(-1)^{F.RL}I_{[t_0,\xi_{1}]}.^{F.RL}I_{[t_0,t_2]}^{\alpha}(^{FC}D^{3\alpha}_{\ast}f_{i.gH})(t_4)\\
\nonumber &\oplus
&~^{F.RL}I_{[\xi_{1},t_1]}.^{F.RL}I_{[\xi_{1},t_3]}^{\alpha}(^{FC}D^{3\alpha}_{\ast}f_{ii.gH})(t_5).
\end{eqnarray}
Finally, the first $(F.RL)$-fractional integral on the right side of
the Eq. $(\ref{eq:4.8})$ obtains as follows
\begin{eqnarray}
\nonumber ^{F.RL}I_{[t_0,\xi]}^{\alpha}(^{FC}D^{\alpha}_{\ast}f_{ii.gH})(t_1)&=&^{FC}D^{\alpha}_{\ast}f_{ii.gH}(t_0)\odot \frac{(\xi-t_0)^{\alpha}}{\Gamma(\alpha+1)}\ominus ~^{FC}D^{2\alpha}_{\ast}f_{i.gH}(t_0)\odot \frac{(t_0-\xi_{1})^{\alpha}}{\Gamma(\alpha+1)}\\
\nonumber &\odot &\frac{(\xi -t_0)^{\alpha}}{\Gamma(\alpha+1)}\oplus ~^{FC}D^{2\alpha}_{\ast}f_{ii.gH}(\xi_{1})\odot[\frac{(\xi -\xi_{1})^{2\alpha}}{\Gamma(2\alpha+1)}-\frac{(t_0-\xi_{1})^{2\alpha}}{\Gamma(2\alpha+1)}]\\
\nonumber &\ominus &(-1)^{F.RL}I_{[t_0,\xi]}^{\alpha}.^{F.RL}I_{[t_0,\xi_{1}]}.^{F.RL}I_{[t_0,t_2]}^{\alpha}(^{FC}D^{3\alpha}_{\ast}f_{i.gH})(t_4)\\
\label{eq:4.14} &\oplus
&^{F.RL}I_{[t_0,\xi]}^{\alpha}.^{F.RL}I_{[\xi_{1},t_1]}.^{F.RL}I_{[\xi_{1},t_3]}^{\alpha}(^{FC}D^{3\alpha}_{\ast}f_{ii.gH})(t_5).
\end{eqnarray}
The only point remaining concerns the behaviour of the second $(F.RL)$-fractional integral on the right side of the Eq. $(\ref{eq:4.8})$. We can now proceed analogously to the first $(F.RL)$-fractional integral:\\
By noting the hypothesis of theorem,
$^{FC}D^{i\alpha}_{\ast}f_{i.gH},~i=2,3$ are
$^{FC}[(i)-gH]$-differentiable on $[\xi,b]$, and the type of
differentiability does not change. By Theorem $\ref{t:3.6}$ we
deduce that
\begin{equation}\label{eq:4.15}
^{FC}D^{\alpha}_{\ast}f_{i.gH}(s_1)=~^{FC}D^{\alpha}_{\ast}f_{i.gH}(\xi)\oplus
~^{F.RL}I_{[\xi
,s_1]}^{\alpha}(^{FC}D^{2\alpha}_{\ast}f_{i.gH})(s_2),
\end{equation}
and
\begin{eqnarray}
\nonumber &&^{FC}D^{2\alpha}_{\ast}f_{i.gH}(s_2)=~^{FC}D^{2\alpha}_{\ast}f_{i.gH}(\xi)\oplus~^{F.RL}I_{[\xi ,s_2]}^{\alpha}(^{FC}D^{3\alpha}_{\ast}f_{i.gH})(s_3).\\
\nonumber &\Rightarrow &^{F.RL}I_{[\xi ,s_1]}^{\alpha}(^{FC}D^{2\alpha}_{\ast}f_{i.gH})(s_2)\\
\label{eq:4.16}&=&^{FC}D^{2\alpha}_{\ast}f_{i.gH}(\xi)\odot
\frac{(s_1-\xi)^{\alpha}}{\Gamma(\alpha+1)}\oplus ~^{F.RL}I_{[\xi
,s_1]}.^{F.RL}I_{[\xi
,s_2]}^{\alpha}(^{FC}D^{3\alpha}_{\ast}f_{i.gH})(s_3).
\end{eqnarray}
Substituting $(\ref{eq:4.16})$ into $(\ref{eq:4.15})$ we obtain
\[^{FC}D^{\alpha}_{\ast}f_{i.gH}(s_1)=~^{FC}D^{\alpha}_{\ast}f_{i.gH}(\xi)\oplus~^{FC}D^{2\alpha}_{\ast}f_{i.gH}(\xi)\odot \frac{(s_1-\xi)^{\alpha}}{\Gamma(\alpha+1)}\oplus~^{F.RL}I_{[\xi ,s_1]}.^{F.RL}I_{[\xi ,s_2]}^{\alpha}(^{FC}D^{3\alpha}_{\ast}f_{i.gH})(s_3).\]
Thus, the second $(F.RL)$-fractional integral on the right side of
the Eq. $(\ref{eq:4.8})$ is as following
\begin{eqnarray}
\nonumber ^{F.RL}I_{[\xi
,t]}^{\alpha}(^{FC}D^{\alpha}_{\ast}f_{i.gH})(s_1)&=&^{FC}D^{\alpha}_{\ast}f_{i.gH
}(\xi)\odot\frac{(t-\xi)^{\alpha}}{\Gamma(\alpha+1)}\oplus~^{FC}D^{2\alpha}_{\ast}f_{i.gH}(\xi)\odot \frac{(t-\xi)^{2\alpha}}{\Gamma(2\alpha+1)}\\
\label{eq:4.17} &\oplus &~^{F.RL}I_{[\xi,t]}.^{F.RL}I_{[\xi
,s_1]}^{\alpha}.^{F.RL}I_{[\xi
,s_2]}^{\alpha}(^{FC}D^{3\alpha}_{\ast}f_{i.gH})(s_3).
\end{eqnarray}
Having disposed of this preliminary step, we can now return to the Eq. $(\ref{eq:4.8})$.\\
By substituting Eq. $(\ref{eq:4.14})$ and Eq. $(\ref{eq:4.17})$, in
Eq. $(\ref{eq:4.8})$, the desired result is achieved.
\end{theor}
\section{Fuzzy Generalized Euler's method}
In this section, we will touch only a few aspects of the fuzzy
generalized Taylor theorem and restrict the discussion to the fuzzy
generalized Euler's method. This case is important enough to be
stated separately. We consider, the following fuzzy fractional
initial value problem
\begin{equation}\label{eq:5.18}
\left\{\begin{array}{ll}
^{FC}D^{\alpha}_{\ast}y_{gH}(t)=f(t,y(t)),&t\in[0,T],
  \\
y(0)=y_0\in \mathbb{R}_{\mathcal{F}},
\end{array}
\right.
\end{equation}
where $f:[0,T]\times\Rb_{\mathcal{F}}\rightarrow\Rb_{\mathcal{F}}$
is continuous and $y(t)$ is an unknown fuzzy function of crisp
variable $t$. Furthermore, $^{FC}D^{\alpha}_{\ast}y_{gH}(t)$ is the
fuzzy fractional derivative $y(t)$ in the Caputo sense of order
$0<\alpha\leq 1$, with the finite set of switching points. Now, by
dividing the interval $[0,T]$ with the step length of $h$, we have
the
partition $\widehat{I}_N=\{0=t_0<t_1<...<t_N=T\}$ where $t_k=kh$ for $k=0,1,2,...,N.$\\
\textbf{Case~I.} Unless otherwise stated we assume that the unique
solution of the fuzzy fractional initial value problem
$(\ref{eq:5.18})$,
$^{FC}D^{2\alpha}_{\ast}y(t)\in\mathcal{C}_{f}([0,T],
\mathbb{R}_{\mathcal{F}} )\cap\mathcal{L}^{\zeta}([0,T],
\mathbb{R}_{\mathcal{F}} )$ is $^{FC}[(i)-gH]$-differentiable such
that the type of differentiability does not change on $[0,T]$.
Consider the fractional Taylor series expansion of the unknown fuzzy
function $y(t)$ about $t_k$, for each $k=0,1,..., N$.
\[y(t_{k+1})=y(t_k)\oplus\frac{(t_{k+1}-t_k)^{\alpha}}{\Gamma(\alpha+1)}\odot~^{FC}D^{\alpha}_{\ast}y_{i.gH}(t_k)\oplus\frac{(t_{k+1}-t_k)^{2\alpha}}{\Gamma(2\alpha+1)}\odot~^{FC}D^{2\alpha}_{\ast}y_{i.gH}(\eta_t),\]
for some points $k$ lie between $t_k$ and $t_{k+1}$. Since $h =
t_{k+1}-t_{k}$, we have
\[y(t_{k+1})=y(t_k)\oplus\frac{h^{\alpha}}{\Gamma(\alpha+1)}\odot~^{FC}D^{\alpha}_{\ast}y_{i.gH}(t_k)\oplus\frac{h^{2\alpha}}{\Gamma(2\alpha+1)}\odot~^{FC}D^{2\alpha}_{\ast}y_{i.gH}(\eta_t),\]
and, y(t) satisfies in problem (5.1), so
\[y(t_{k+1})=y(t_k)\oplus\frac{h^{\alpha}}{\Gamma(\alpha+1)}\odot~f(t_k,y(t_k))\oplus\frac{h^{2\alpha}}{\Gamma(2\alpha+1)}\odot~^{FC}D^{2\alpha}_{\ast}y_{i.gH}(\eta_t),\]
\begin{eqnarray}
\nonumber &&\mathcal{H}(y(t_{k+1}),y(t_k)\oplus\frac{h^{\alpha}}{\Gamma(\alpha+1)}\odot~f(t_k,y(t_k))\oplus\frac{h^{2\alpha}}{\Gamma(2\alpha+1)}\odot~^{FC}D^{2\alpha}_{\ast}y_{i.gH}(\eta_t))\\
\nonumber &\leq
&\mathcal{H}(y(t_{k+1}),y(t_k)\oplus\frac{h^{\alpha}}{\Gamma(\alpha+1)}\odot~f(t_k,y(t_k))+\mathcal{H}(0,\frac{h^{2\alpha}}{\Gamma(2\alpha+1)}\odot~^{FC}D^{2\alpha}_{\ast}y_{i.gH}(\eta_t)),
\end{eqnarray}
as $h\rightarrow 0$ since
\begin{eqnarray}
\nonumber &&\mathcal{H}(y(t_{k+1}),y(t_k)\oplus\frac{h^{\alpha}}{\Gamma(\alpha+1)}\odot f(t_k,y(t_k))\rightarrow 0,\\
\nonumber
&&\mathcal{H}(0,\frac{h^{2\alpha}}{\Gamma(2\alpha+1)}\odot~^{FC}D^{2\alpha}_{\ast}y_{i.gH}(\eta_t))\rightarrow
0,
\end{eqnarray}
we conclude that
\begin{eqnarray}
\nonumber &&\mathcal{H}(y(t_{k+1}),y(t_k)\oplus\frac{h^{\alpha}}{\Gamma(\alpha+1)}\odot f(t_k,y(t_k))+\mathcal{H}(0,\frac{h^{2\alpha}}{\Gamma(2\alpha+1)}\odot~^{FC}D^{2\alpha}_{\ast}y_{i.gH}(\eta_t))\rightarrow 0,\\
\nonumber &\Rightarrow
&\mathcal{H}_d(y(t_{k+1}),y(t_k)\oplus\frac{h^{\alpha}}{\Gamma(\alpha+1)}\odot~f(t_k,y(t_k))\oplus\frac{h^{2\alpha}}{\Gamma(2\alpha+1)}\odot~^{FC}D^{2\alpha}_{\ast}y_{i.gH}(\eta_t))\rightarrow
0.
\end{eqnarray}
Thus, for sufficiently small $h$ we find that
\[y(t_{k+1})\approx y(t_k)\oplus\frac{h^{\alpha}}{\Gamma(\alpha+1)}\odot~f(t_k,y(t_k)),\]
and finally we get
\begin{equation}\label{eq:5.19}
\left\{\begin{array}{ll}
  y_0=y_0,
  \\
  y_{k+1}=y_k\oplus\frac{h^{\alpha}}{\Gamma(\alpha+1)}\odot f(t_k,y_k),&k=0,1, ...,N-1.
\end{array}
\right.
\end{equation}
\textbf{Case~II.} Assume that
$^{FC}D^{2\alpha}_{\ast}y(t)\in\mathcal{C}_f([0,T],\mathbb{R}_{\mathcal{F}})$
is $^{FC}[(ii)-gH]$-differentiable such that the type of
differentiability does not change on $[0, T]$. So the fractional
Taylor's series expansion of $y(t)$ about the point $t_k$ at
$t_{k+1}$ is
\[y(t_{k+1})=y(t_k)\ominus(-1)\frac{(t_{k+1}-t_k)^{\alpha}}{\Gamma(\alpha+1)}\odot~^{FC}D^{\alpha}_{\ast}y_{ii.gH}(t_k)\ominus(-1)\frac{(t_{k+1}-t_k)^{2\alpha}}{\Gamma(2\alpha+1)}\odot~^{FC}D^{2\alpha}_{\ast}y_{ii.gH}(\eta_t).\]
According to the process described in Case \textbf{I}, the
generalized Euler's method takes the form
\begin{equation}\label{eq:5.20}
\left\{\begin{array}{ll}
  y_0=y_0,
  \\
  y_{k+1}=y_k\ominus(-1)\frac{h^{\alpha}}{\Gamma(\alpha+1)}\odot f(t_k,y_k),&k=0,1, ...,N-1.
\end{array}
\right.
\end{equation}
\textbf{Case~III.} Let us suppose that $t_0=0, t_1, ..., t_j, \zeta,
t_{j+1}, ..., t_N=T$ is a partition of interval $[0,T]$ and $y(t)$
has a switching point in $\zeta\in[0,T]$ of type \textbf{I}. So
according to Eqs.  $(\ref{eq:5.19})$ and $(\ref{eq:5.20})$, we have
\begin{equation}\label{eq:5.21}
\left\{\begin{array}{ll}
  y_0=y_0,
  \\
  y_{k+1}=y_k\oplus\frac{h^{\alpha}}{\Gamma(\alpha+1)}\odot f(t_k,y_k),&k=0,1, ...,j.
  \\
  y_{k+1}=y_k\ominus(-1)\frac{h^{\alpha}}{\Gamma(\alpha+1)}\odot f(t_k,y_k),&k=j+1,j+2, ...,N-1.
\end{array}
\right.
\end{equation}
\textbf{Case~IV.} Consider $y(t)$ has a switching point type
\textbf{II} in $\zeta\in[0,T]$ such that $t_0, t_1, ..., t_j, \zeta,
t_{j+1}, ..., t_N$ is a partition of interval $[0,T]$. Hence by Eqs.
$(\ref{eq:5.19})$ and $(\ref{eq:5.20})$, we conclude that
\begin{equation}\label{eq:5.22}
\left\{\begin{array}{ll}
   y_0=y_0,
  \\
 y_{k+1}=y_k\ominus(-1)\frac{h^{\alpha}}{\Gamma(\alpha+1)}\odot f(t_k,y_k),&k=0,1, ...,j.
  \\
   y_{k+1}=y_k\oplus\frac{h^{\alpha}}{\Gamma(\alpha+1)}\odot f(t_k,y_k),&k=j+1,j+2, ...,N-1.
\end{array}
\right.
\end{equation}
Our next concern will be the behavior of the fuzzy generalized Euler
method.
\section{Analysis of the Fuzzy Generalized Euler's method}\label{sec5}
In this section, the local and the global truncation errors of the
fuzzy generalized Euler's method are illustrated. So by applying
them the consistence, the convergence and the stability of the
presented method are proved. Also, several definitions and concepts
of the fuzzy generalized Euler's method are presented under
$^{FC}[gH]$-differentiability \cite{an}.
\subsection{Local Truncation Error, Consistent}
Consider the unique solution of the fuzzy fractional initial value
problem $(\ref{eq:5.18}):$
\begin{definition}\label{d:6.1}
If $y(t)$ is $^{FC}[(i)-gH]$-differentiable on $[0,T]$ and the type
of differentiability does not change, now we define the residual
$\mathcal{R}_k$  as
\[\mathcal{R}_k=y(t_{k+1})\ominus _{gH}\left(y(t_k)\oplus\frac{h^{\alpha}}{\Gamma(\alpha+1)}\odot f(t_k,y(t_k))\right),\]
and if $y(t)$ is $^{FC}[(ii)-gH]$-differentiable on $[0,T]$, we have
\[\mathcal{R}_k=y(t_{k+1})\ominus _{gH}\left(y(t_k)\ominus(-1)\frac{h^{\alpha}}{\Gamma(\alpha+1)}\odot f(t_k,y(t_k))\right).\]
On the other hand, the local truncation error (L. T. E.) $(\tau_k)$
is defined as
\[\tau_k=\frac{1}{h}\mathcal{R}_k,\]
and the fuzzy generalized Euler's method is said to be consistent if
\[\lim_{h\rightarrow0}\max_{t_k\leq T} \mathcal{H}(\tau_k,0)=0.\]
Therefore, due to the type of differentiability of $y(t)$ for
$\eta_k\in[t_k,t_{k+1}]$, the
residual $(\mathcal{R}_k)$ and the L. T. E. $(\tau_k)$ are defined as follows:\\

$\bullet~~\begin{array}{ll}
&~~~\mathcal{R}_k=\frac{h^{2\alpha}}{\Gamma(2\alpha+1)}\odot
~^{FC}D^{2\alpha}_{\ast}y_{i.gH}(\eta_k),\\
^{FC}[(i)-gH]-differentiability\Rightarrow
\\
&~~~\tau_k=\frac{h^{2\alpha-1}}{\Gamma(2\alpha+1)}\odot
~^{FC}D^{2\alpha}_{\ast}y_{i.gH}(\eta_k),
\end{array}$
\\\\

$\bullet~~\begin{array}{ll}
&~~~\mathcal{R}_k=\ominus(-1)\frac{h^{2\alpha}}{\Gamma(2\alpha+1)}\odot
~^{FC}D^{2\alpha}_{\ast}y_{ii.gH}(\eta_k),\\^{FC}[(ii)-gH]-differentiability\Rightarrow
\\
&~~~\tau_k=\ominus(-1)\frac{h^{2\alpha-1}}{\Gamma(2\alpha+1)}\odot
~^{FC}D^{2\alpha}_{\ast}y_{ii.gH}(\eta_k).
\end{array}$\\
\end{definition}
$\triangleleft$ Investigating the \textbf{consistence} of the fuzzy generalized Euler's method\textbf{:}\\
For this purpose, assume that
$\mathcal{H}(^{FC}D^{2\alpha}_{\ast}y_{ii.gH}(\eta_k),0)\leq M$. We
have two following steps:

\textbf{step~I.} If $y(t)$ be $^{FC}[(i)-gH]$-differentiable, then
\begin{eqnarray}
\nonumber \lim_{h\rightarrow 0}\max_{t_k\leq T} \mathcal{H}(\tau_k,0)&=&\lim_{h\rightarrow 0}\max_{t_k\leq T} \mathcal{H}(\frac{h^{2\alpha-1}}{\Gamma(2\alpha+1)}\odot ~^{FC}D^{2\alpha}_{\ast}y_{i.gH}(\eta_k),0)\\
\nonumber &=&\lim_{h\rightarrow 0}\frac{h^{2\alpha-1}}{\Gamma(2\alpha+1)}\max_{t_k\leq T}\mathcal{H}(^{FC}D^{2\alpha}_{\ast}y_{i.gH}(\eta_k),0)\\
\nonumber &\leq &\lim_{h\rightarrow 0}\frac{h^{2\alpha-1}}{\Gamma(2\alpha+1)}.M=0.
\end{eqnarray}

\textbf{step II.} The same conclusion can be drawn for the
$^{FC}[(ii)-gH]$-differentiability of $y(t)$, so
\begin{eqnarray}
\nonumber\lim_{h\rightarrow 0}\max_{t_k\leq T} \mathcal{H}(\tau_k,0)&=&\lim_{h\rightarrow 0}\max_{t_k\leq T} \mathcal{H}(\ominus(-1)\frac{h^{2\alpha-1}}{\Gamma(2\alpha+1)}\odot ~^{FC}D^{2\alpha}_{\ast}y_{ii.gH}(\eta_k),0)\\
\nonumber &=&\lim_{h\rightarrow 0}\mid(-1)\frac{h^{2\alpha-1}}{\Gamma(2\alpha+1)}\mid\max_{t_k\leq T} \mathcal{H}(\ominus ~^{FC}D^{2\alpha}_{\ast}y_{ii.gH}(\eta_k),0)\\
\nonumber &=&\lim_{h\rightarrow
0}\frac{h^{2\alpha-1}}{\Gamma(2\alpha+1)}.\mathcal{H}(^{FC}D^{2\alpha}_{\ast}y_{ii.gH}(\eta_k),0)\leq\lim_{h\rightarrow
0}\frac{h^{2\alpha-1}}{\Gamma(2\alpha+1)}.M=0.
\end{eqnarray}
Thus, note that we have actually proved that the fuzzy generalized
Euler's method is consistent as long as the solution belongs to
$\mathcal{C}_f([0,T],~\Rb_{\mathcal{F}}).$
\subsection{Global Truncation Error, Convergence}
\begin{lemma}\label{l:6.1}\cite{e}
$\forall\emph{z}\in\Rb,~1+\emph{z}\leq e^{\emph{z}}$.
\end{lemma}
\begin{definition}\cite{g}\label{d:6.2}
The global truncation error is the agglomeration of the local
truncation error over all the iterations, assuming perfect knowledge
of the true solution at the initial time step.

In the fuzzy fractional initial value problem $(\ref{eq:5.18})$,
assume that $y(t)$ is $^{FC}[(i)-gH]$-differentiable, then the
global truncation error is
\begin{eqnarray}
\nonumber e_{k+1}&=&y(t_{k+1})\ominus _{gH}y_{k+1}=y(t_{k+1})\ominus _{gH}[y_0\oplus\frac{h^{\alpha}}{\Gamma(\alpha+1)}\odot f(t_0,y_0)\\
\nonumber &\oplus &\frac{h^{\alpha}}{\Gamma(\alpha+1)}\odot
f(t_1,y_1)\oplus...\oplus\frac{h^{\alpha}}{\Gamma(\alpha+1)}\odot
f(t_k,y_k))],
\end{eqnarray}
and for the $^{FC}[(ii)-gH]$-differentiability of $y(t)$, we have
\begin{eqnarray}
\nonumber e_{k+1}&=&y(t_{k+1})\ominus _{gH}y_{k+1}=y(t_{k+1})\ominus _{gH}[y_0\ominus(-1)\frac{h^{\alpha}}{\Gamma(\alpha+1)}\odot f(t_0,y_0)\\
\nonumber &\ominus &(-1)\frac{h^{\alpha}}{\Gamma(\alpha+1)}\odot
f(t_1,y_1)\ominus(-1)...\ominus(-1)\frac{h^{\alpha}}{\Gamma(\alpha+1)}\odot
f(t_k,y_k))].
\end{eqnarray}
\end{definition}
\begin{definition}\label{d:6.3}
If global truncation error leads to zero as the step size goes to
zero, the numerical method is convergent, i.e.
\[\nonumber \lim_{h\rightarrow0}\max_{k} \mathcal{H}(e_{k+1},0)=0,~\Rightarrow ~\lim_{h\rightarrow0}\max_{k} \mathcal{H}\left(y(t_{k+1}),y_{k+1}\right)=0.\]
In this case, the numerical solution converges to the exact
solution.
\end{definition}
$\triangleleft$ Investigating the \textbf{convergence} of the fuzzy generalized Euler's method\textbf{:}\\
To suppose that $^{FC}D^{2\alpha}_{\ast}y(t)$ exists and $f(t,y)$
satisfies in Lipschitz condition on the $\{(t,y)\mid t\in[0,
p],~y\in \overline{B}(y_0,q),~p,~q>0\}$, the research on this
subject will be divided into two steps:

\textbf{step I.} Suppose that $y(t)$ is
$^{FC}[(i)-gH]$-differentiable, now by using Eq. $(\ref{eq:5.19})$
and assumption
$r_k=\frac{h^{2\alpha}}{\Gamma(2\alpha+1)}\odot~^{FC}D^{2\alpha}_{\ast}y_{i.gH}(t_k)$
the exact solution of the fuzzy fractional initial value problem
$(\ref{eq:5.18})$ satisfies
\[y(t_{k+1})=y(t_k)\oplus\frac{h^{\alpha}}{\Gamma(\alpha+1)}\odot~f(t_k,y(t_k))\oplus~r_k.\]
Subtracting the above equation from Eq. $(\ref{eq:5.19})$, deduces
\[\mathcal{H}\left(y(t_{k+1}),y_{k+1}\right)=\mathcal{H}\left(y(t_{k}),y_{k}\right)+\frac{h^{\alpha}}{\Gamma(\alpha+1)}.\mathcal{H}\left(f(t_k,y(t_k)),f(t_k,y_k)\right)+\mathcal{H}(r_k,0).\]
Since $f(t,y(t))$ satisfies in Lipschitz condition
\[\mathcal{H}\left(f(t_k,y(t_k)),f(t_k,y_k)\right)\leq\ell_k\mathcal{H}_d\left(y(t_{k}),y_{k}\right),\]
and this inequality obey,
\begin{equation}\label{eq:6.23}
\mathcal{H}\left(y(t_{k+1}),y_{k+1}\right)\leq(1+\frac{h^{\alpha}}{\Gamma(\alpha+1)}.\ell_k)\mathcal{H}\left(y(t_{k}),y_{k}\right)+\mathcal{H}(r_k,0).
\end{equation}
From now on let
\[\ell=\max_{0\leq k\leq N}~\ell_k,~~~~~~~r=\max_{0\leq k\leq N}~\mathcal{H}(r_k,0),\]
thus, the Eq. $(\ref{eq:6.23})$ can be written as
\[\mathcal{H}\left(y(t_{k+1}),y_{k+1}\right)\leq(1+\frac{h^{\alpha}}{\Gamma(\alpha+1)}.\ell)\mathcal{H}\left(y(t_{k}),y_{k}\right)+r.\]
Since the inequality holds for all $k$, it follows that
\begin{eqnarray}
\nonumber \mathcal{H}\left(y(t_{k+1}),y_{k+1}\right)&\leq &(1+\frac{h^{\alpha}}{\Gamma(\alpha+1)}.\ell)\left[(1+\frac{h^{\alpha}}{\Gamma(\alpha+1)}.\ell)\mathcal{H}\left(y(t_{k-1}),y_{k-1}\right)+r\right]+r\\
\nonumber
&=&(1+\frac{h^{\alpha}}{\Gamma(\alpha+1)}.\ell)^2.\mathcal{H}\left(y(t_{k-1}),y_{k-1}\right)+r\left[1+(1+\frac{h^{\alpha}}{\Gamma(\alpha+1)}.\ell)\right].
\end{eqnarray}
In the same trend, finds
\begin{eqnarray}
\nonumber &&\mathcal{H}\left(y(t_{k+1}),y_{k+1}\right)\\
\nonumber &\leq
&(1+\frac{h^{\alpha}}{\Gamma(\alpha+1)}.\ell)^{k+1}.\mathcal{H}\left(y(t_{0}),y_{0}\right)+r\left[1+(1+\frac{h^{\alpha}}{\Gamma(\alpha+1)}.\ell)+...+(1+\frac{h^{\alpha}}{\Gamma(\alpha+1)}.\ell)^k\right].
\end{eqnarray}
Using the formula for the sum of a geometric series, we obtain
\[\sum_{i=0}^k(1+\frac{h^{\alpha}}{\Gamma(\alpha+1)}.\ell)^i=\frac{(1+\frac{h^{\alpha}}{\Gamma(\alpha+1)}.\ell)^{k+1}-1}{\frac{h^{\alpha}}{\Gamma(\alpha+1)}.\ell},\]
which leads to the following inequality
\[\mathcal{H}\left(y(t_{k+1}),y_{k+1}\right)\leq (1+\frac{h^{\alpha}}{\Gamma(\alpha+1)}.\ell)^{k+1}\mathcal{H}
\left(y(t_{0}),y_{0}\right)+\frac{r\Gamma(\alpha+1)}{h^{\alpha}\ell}\left[(1+\frac{h^{\alpha}}{\Gamma(\alpha+1)}.\ell)^{k+1}-1\right].\]
On the other hand, Lemma $\ref{l:6.1}$ and
$0\leq(k+1)h^{\alpha}\leq T$ (for
$(k+1)\leq (N-1)$) yields
\[\mathcal{H}\left(y(t_{k+1}),y_{k+1}\right)\leq e^{\frac{\ell T}{\Gamma(\alpha+1)}}\mathcal{H}\left(y(t_{0}),y_{0}\right)+\frac{r\Gamma(\alpha+1)}{h^{\alpha}\ell}[e^{\frac{\ell T}{\Gamma(\alpha+1)}}-1],\]
given that $\mathcal{H}\left(y(t_{0}),y_{0}\right)=0$ and
\[r=\max_{0\leq k\leq N-1}\mathcal{H}(r_k,0)=\frac{h^{2\alpha}}{\Gamma(2\alpha+1)}\max_{0\leq t\leq T}\mathcal{H}( ^{FC}D^{2\alpha}_{\ast}y_{i.gH}(t),0),\]
immediately concludes that
\[\mathcal{H}\left(y(t_{k+1}),y_{k+1}\right)\leq\frac{h^{\alpha}\Gamma(\alpha+1)}{\ell\Gamma(2\alpha+1)}[e^{\frac{\ell T}{\Gamma(\alpha+1)}}-1]\max_{0\leq t\leq T}\mathcal{H}( ^{FC}D^{2\alpha}_{\ast}y_{i.gH}(t),0).\]
So, $\lim_{h\rightarrow
0}\mathcal{H}\left(y(t_{k+1}),y_{k+1}\right)\rightarrow 0$ and in
this step, the fuzzy generalized Euler's method is convergent.

\textbf{step II.} To estimate the step II, consider $y(t)$ is
$^{FC}[(ii)-gH]$-differentiable, by using Eq. (\ref{eq:5.20}) and
let

$r_k=\ominus(-1)\frac{h^{2\alpha}}{\Gamma(2\alpha+1)}\odot~^{FC}D^{2\alpha}_{\ast}y_{ii.gH}(t_k)$,
the exact solution of the Eq. (\ref{eq:5.18}) satisfies
\[y(t_{k+1})=y(t_k)\ominus(-1)\frac{h^{\alpha}}{\Gamma(\alpha+1)}\odot~f(t_k,y(t_k))\oplus~r_k.\]
\[\Rightarrow \mathcal{H}\left(y(t_{k+1}),y_{k+1}\right)=\mathcal{H}\left(y(t_{k}),y_{k}\right)+\frac{h^{\alpha}}{\Gamma(\alpha+1)}.[\mathcal{H}\left(f(t_k,y_k)\ominus_{gH}f(t_k,y(t_k)),0\right)]+\mathcal{H}(r_k,0).\]
The inequality
\[\mathcal{H}\left(f(t_k,y_k)\ominus_{gH}f(t_k,y(t_k)),0\right)=\mathcal{H}\left(f(t_k,y_k),f(t_k,y(t_k))\right)
\leq\ell_k\mathcal{H}\left(y(t_{k}),y_{k}\right)\],  which is the
conclusion of Lipschitz condition, implies that
\begin{equation}\label{eq:6.24}
\mathcal{H}\left(y(t_{k+1}),y_{k+1}\right)\leq(1-\frac{h^{\alpha}}{\Gamma(\alpha+1)}.\ell_k)\mathcal{H}\left(y(t_{k}),y_{k}\right)+\mathcal{H}(r_k,0).
\end{equation}
Now, assume that
\[\ell=\max_{0\leq k\leq N-1}~\ell_k,~~~~~r=\max_{0\leq k\leq N-1}~\mathcal{H}(r_k,0),\]
 and rewrite Eq. $(\ref{eq:6.24})$ as
\[\mathcal{H}\left(y(t_{k+1}),y_{k+1}\right)\leq(1-\frac{h^{\alpha}}{\Gamma(\alpha+1)}.\ell)\mathcal{H}\left(y(t_{k}),y_{k}\right)+r.\]
Since the inequality holds for all $k$, we get
\begin{eqnarray}
\nonumber \mathcal{H}\left(y(t_{k+1}),y_{k+1}\right)&\leq &(1-\frac{h^{\alpha}}{\Gamma(\alpha+1)}.\ell)\left[(1-\frac{h^{\alpha}}{\Gamma(\alpha+1)}.\ell)\mathcal{H}\left(y(t_{k-1}),y_{k-1}\right)+r\right]+r\\
\nonumber
&=&(1-\frac{h^{\alpha}}{\Gamma(\alpha+1)}.\ell)^2.\mathcal{H}\left(y(t_{k-1}),y_{k-1}\right)+r\left[1+(1-\frac{h^{\alpha}}{\Gamma(\alpha+1)}.\ell)\right].
\end{eqnarray}
Repeated application of the above inequality enables us to write
\begin{eqnarray}
\nonumber &&\mathcal{H}\left(y(t_{k+1}),y_{k+1}\right)\\
\nonumber &\leq
&(1-\frac{h^{\alpha}}{\Gamma(\alpha+1)}.\ell)^{k+1}.\mathcal{H}\left(y(t_{0}),y_{0}\right)+r\left[1+(1-\frac{h^{\alpha}}{\Gamma(\alpha+1)}.\ell)+...+(1-\frac{h^{\alpha}}{\Gamma(\alpha+1)}.\ell)^k\right].
\end{eqnarray}
Obviously, this sum is a geometric series, so we have
\[\sum_{i=0}^k(1-\frac{h^{\alpha}}{\Gamma(\alpha+1)}.\ell)^i=\frac{1-(1-\frac{h^{\alpha}}{\Gamma(\alpha+1)}.\ell)^{k+1}}{\frac{h^{\alpha}}{\Gamma(\alpha+1)}.\ell},\]
that resulted to
\begin{equation}\label{eq:6.25}
\mathcal{H}\left(y(t_{k+1}),y_{k+1}\right)\leq
(1-\frac{h^{\alpha}}{\Gamma(\alpha+1)}.\ell)^{k+1}\mathcal{H}\left(y(t_{0}),y_{0}\right)+\frac{r\Gamma(\alpha+1)}{h^{\alpha}\ell}\left[1-(1-\frac{h^{\alpha}}{\Gamma(\alpha+1)}.\ell)^{k+1}\right].
\end{equation}
With $\emph{z}=-\frac{h^{\alpha}}{\Gamma(\alpha+1)}.\ell$ in Lemma
$\ref{l:6.1}$ concludes that
\[(1-\frac{h^{\alpha}}{\Gamma(\alpha+1)}.\ell)^{k+1}\leq e^{-\frac{h^{\alpha}}{\Gamma(\alpha+1)}.\ell(k+1)}\leq e^{-\frac{\ell T}{\Gamma(\alpha+1)}},\]
where $0\leq(k+1)h^{\alpha}\leq T$ for
$(k+1)\leq (N-1).$ Thus in Eq. $(\ref{eq:6.25})$, we obtain
\[\mathcal{H}\left(y(t_{k+1}),y_{k+1}\right)\leq e^{-\frac{\ell T}{\Gamma(\alpha+1)}}\mathcal{H}\left(y(t_{0}),y_{0}\right)+\frac{r\Gamma(\alpha+1)}{h^{\alpha}\ell}[1-e^{-\frac{\ell T}{\Gamma(\alpha+1)}}].\]
Moreover
\[r=\max_{0\leq k\leq N-1}\mathcal{H}(r_k,0)=-\frac{h^{2\alpha}}{\Gamma(2\alpha+1)}\max_{0\leq t\leq T}\mathcal{H}( ^{FC}D^{2\alpha}_{\ast}y_{ii.gH}(t),0),\]
and the accuracy of the initial value, concludes that
$\mathcal{H}\left(y(t_{0}),y_{0}\right)=0$, so
\[\mathcal{H}\left(y(t_{k+1}),y_{k+1}\right)\leq -\frac{h^{\alpha}\Gamma(\alpha+1)}{\ell\Gamma(2\alpha+1)}[1-e^{-\frac{\ell T}{\Gamma(\alpha+1)}}]\max_{0\leq t\leq T}\mathcal{H}( ^{FC}D^{2\alpha}_{\ast}y_{ii.gH}(t),0).\]
Now, letting ${h\rightarrow 0}$ then
$\mathcal{H}\left(y(t_{k+1}),y_{k+1}\right)\rightarrow 0$ which is
the desired conclusion and we can say that  the fuzzy generalized
Euler's method is convergent.
\subsection{Stability}
Now, the stability of the presented method is illustrated. For this
aim, the following definition is presented as
\begin{definition}\label{d:6.4}
Assume that $y_{k+1},~k+1\geq 0$ is the solution of fuzzy
generalized Euler's method where $y_0\in \mathbb{R}_{\mathcal{F}}$
and also $z_{k+1}$ is the solution of the same numerical method
where $z_0 = y_0 \oplus \delta_{0}\in \mathbb{R}_{\mathcal{F}}$
shows its perturbed fuzzy initial condition. The fuzzy generalized
Euler's method is stable if there exists positive constant
$\widehat{h}$ and $\mathcal{K}$ such that
\[\forall~(k+1)h^{\alpha}\leq T,~k+1< N-1,~h\in (0,\widehat{h})\Rightarrow\mathcal{H}\left(z_{k+1},y_{k+1}\right)\leq\mathcal{K}\delta\]
whenever $\mathcal{H}\left(\delta_0,0\right)\leq\delta.$
\end{definition}
$\triangleleft$ Investigating the \textbf{stability} of the fuzzy generalized Euler's method\textbf{:}\\
The proof falls naturally into two steps:

\textbf{step I.} If $y(t)$ is $^{FC}[(i)-gH]$-differentiable, by
using Eq. $(\ref{eq:5.19})$ the perturbed problem is in the
following form
\begin{equation}\label{eq:6.26}
z_{k+1}=z_k\oplus \frac{h^{\alpha}}{\Gamma(\alpha+1)}\odot f(t_k,
z_k),~~~z_0=y_0\oplus \delta_0.
\end{equation}
So, considering the Eqs. $(\ref{eq:5.19})$ and $(\ref{eq:6.26})$,
gets
\[\mathcal{H}\left(z_{k+1},y_{k+1}\right)\leq\mathcal{H}\left(z_{k},y_{k}\right)+\frac{h^{\alpha}}{\Gamma(\alpha+1)}\mathcal{H}\left(f(t_k,z_{k}),f(t_k,y_{k})\right).\]
By using the properties of Hausdorff metric and the Lipschitz
condition that expressed in Section 2, we have
\[\mathcal{H}\left(z_{k+1},y_{k+1}\right)\leq(1+\frac{h^{\alpha}}{\Gamma(\alpha+1)}\ell)\mathcal{H}\left(z_{k},y_{k}\right).\]
Continuing this process and iterating the inequality, leads to the
following relation
\[\mathcal{H}\left(z_{k+1},y_{k+1}\right)\leq(1+\frac{h^{\alpha}}{\Gamma(\alpha+1)}\ell)^{k+1}\mathcal{H}\left(z_0,y_0\right).\]
Now, the Lemma $\ref{l:6.1}$ implies
\[\mathcal{H}\left(z_{k+1},y_{k+1}\right)\leq e^{\frac{h^{\alpha}}{\Gamma(\alpha+1)}.\ell(k+1)}\mathcal{H}\left(z_0\ominus _{gH}y_0,0\right),\]
and finally
\[\mathcal{H}\left(z_{k+1},y_{k+1}\right)\leq e^{\frac{\ell T}{\Gamma(\alpha+1)}}\mathcal{H}(\delta _0,0)\leq\mathcal{K}\delta,\]
where $\mathcal{K}=e^{\frac{\ell T}{\Gamma(\alpha+1)}}$ and for $k+1<N-1\Rightarrow
h^{\alpha}(k + 1)\leq T$. In this case, it
is obvious the stability of the fuzzy generalized Euler's method.

\textbf{step II.} The same proof obtains when we consider the
assumption $^{FC}[(ii)-gH]$-differentiability of $y(t)$. The
numerical method $(\ref{eq:5.20})$ is applied to perturbation
problem. So we get
\begin{equation}\label{eq:6.27}
z_{k+1}=z_k\ominus(-1)\frac{h^{\alpha}}{\Gamma(\alpha+1)}\odot
f(t_k, z_k),~~~z_0=y_0\oplus \delta_0.
\end{equation}
According to the Eqs. $(\ref{eq:5.20})$ and $(\ref{eq:6.27})$, we
have
\[\mathcal{H}\left(z_{k+1},y_{k+1}\right)\leq\mathcal{H}\left(z_{k},y_{k}\right)-\frac{h^{\alpha}}{\Gamma(\alpha+1)}\mathcal{H}\left(f(t_k,z_{k}),f(t_k,y_{k})\right),\]
which we have been working under the assumption that specifications
of the Hausdorff metric are satisfied. Using the Lipschitz condition
can be concluded that
\[\mathcal{H}\left(z_{k+1},y_{k+1}\right)\leq(1-\frac{h^{\alpha}}{\Gamma(\alpha+1)}\ell)\mathcal{H}\left(z_{k},y_{k}\right).\]
Repeating with the inequality and applying Lemma $\ref{l:6.1}$, lead
us to the following inequality
\begin{eqnarray}
\nonumber \mathcal{H}\left(z_{k+1},y_{k+1}\right)&\leq &(1-\frac{h^{\alpha}}{\Gamma(\alpha+1)}\ell)^{k+1}\mathcal{H}\left(z_0,y_0\right)\\
\nonumber &\leq &e^{-\frac{h^{\alpha}}{\Gamma(\alpha+1)}.\ell(k+1)}\mathcal{H}\left(z_0\ominus _{gH}y_0,0\right)\\
\nonumber &\leq &e^{-\frac{\ell T}{\Gamma(\alpha+1)}}\mathcal{H}(\delta
_0,0)\leq\mathcal{K}\delta,
\end{eqnarray}
where $\mathcal{K}=e^{-\frac{\ell T}{\Gamma(\alpha+1)}}$. For the general case, above
analysis, just amounts to the fact that the fuzzy generalized Euler
method is a stable approach.
\section{Numerical Simulations}
In this section, several examples of the fractional differential
equations are solved by using the full fuzzy generalized Euler
method. Also, the numerical results are demonstrated on some tables
for different values of $h$ and $t$.
\begin{example}\label{ex1} Let us consider the following initial value
problem
$$
^{FC}D^{\alpha}_{\ast}y(t)=(0,1,1.5)\odot\Gamma (\alpha+1),
~~~~~0\leq t\leq 1,
$$
where $y(0)=0$ and $y(t) = (0,1,1.5)\odot t^{\alpha}$ is the exact
$^{FC}[i-gH]$-differentiable solution of problem. In order to find
the numerical results we should construct the following iterative
formula as
$$
y_{k+1} = y_k \oplus \frac{h^{\alpha}}{\Gamma(\alpha+1)}\odot
\left[(0,1,1.5)\odot\Gamma (\alpha+1)\right],~~~~~k=0,1,\cdots, N-1.
$$
In Table \ref{t1}, the numerical results for different values of $t,
\alpha$ and $h$ are demonstrated. In Fig. \ref{f1}, the exact
solution and the Caputo gH-derivative for $\alpha=0.6$ are
demonstrated.
\begin{table}
\caption{Numerical results of Example \ref{ex1} for various $t,
\alpha$ and $h$.}\label{t1}
 \centering
\scalebox{0.6}{
\begin{tabular}{|c|l|l|l|l|l|l|}
  \hline
  \multicolumn{3}{ |c| }{$\alpha =0.3$}& \multicolumn{2}{ |c| }{$\alpha =0.6$ }& \multicolumn{2}{ |c| }{$\alpha =0.9$ } \\
 \hline
$t$   & $h=0.2$ & $h=0.02$ &  $h=0.2$ & $h=0.02$ &   $h=0.2$ & $h=0.02$ \\
\hline
0.1&(0,0.617034,0.925551)&(0,0.309249,0.463874)&(0,0.380731,0.571096)&(0,0.0956352,0.143453)&(0,0.234924,0.352386)&(0,0.0295752,0.0443627)\\
0.2&(0,1.23407,1.8511)&(0,0.618499,0.927748)&(0,0.761462,1.14219)&(0,0.19127,0.286906)&(0,0.469848,0.704771)&(0,0.0591503,0.0887255)\\
0.3&(0,1.8511,2.77665)&(0,0.927748,1.39162)&(0,1.14219,1.71329)&(0,0.286906,0.430359)&(0,0.704771,1.05716)&(0,0.0887255,0.133088)\\
0.4&(0,2.46814,3.7022)&(0,1.237,1.8555)&(0,1.52292,2.28438)&(0,0.382541,0.573811)&(0,0.939695,1.40954)&(0,0.118301,0.177451)\\
0.5&(0,3.08517,4.62775)&(0,1.54625,2.31937)&(0,1.90365,2.85548)&(0,0.478176,0.717264)&(0,1.17462,1.76193)&(0,0.147876,0.221814)\\
0.6&(0,3.7022,5.5533)&(0,1.8555,2.78325)&(0,2.28438,3.42658)&(0,0.573811,0.860717)&(0,1.40954,2.11431)&(0,0.177451,0.266176)\\
0.7&(0,4.31924,6.47886)&(0,2.16475,3.24712)&(0,2.66512,3.99767)&(0,0.669447,1.00417)&(0,1.64447,2.4667)&(0,0.207026,0.310539)\\
0.8&(0,4.93627,7.40441)&(0,2.474,3.71099)&(0,3.04585,4.56877)&(0,0.765082,1.14762)&(0,1.87939,2.81909)&(0,0.236601,0.354902)\\
0.9&(0,5.5533,8.32996)&(0,2.78325,4.17487)&(0,3.42658,5.13987)&(0,0.860717,1.29108)&(0,2.11431,3.17147)&(0,0.266176,0.399265)\\
1.0&(0,6.17034,9.25551)&(0,3.09249,4.63874)&(0,3.80731,5.71096)&(0,0.956352,1.43453)&(0,2.34924,3.52386)&(0,0.295752,0.443627)\\
 \hline
 \end{tabular}
 }
\end{table}
\end{example}
\begin{example}\label{ex2} Consider the following problem
$$
^{FC}D^{\alpha}_{\ast}y(t) =(-1)\odot y(t), ~~~~~0\leq t\leq 1,
$$
where $y(0) = (0,1,2)$ and the exact $^{FC}[ii-gH]$-differentiable
solution of problem is in the form $y(t) = (0,1,2)\odot
E_{\alpha}(-t^{\alpha})$. In order to solve the mentioned problem
the following formula should be applied as
$$
\begin{array}{l}
  y_0=(0,1,2), \\
  y_{k+1}=y_k \ominus_{gH} \frac{h^{\alpha}}{\Gamma(\alpha+1)} \odot y_k,
~~~~~k=0,1,\cdots,N-1.
\end{array}
$$
The numerical results based on the presented method are obtained in
Table \ref{t2} for various $t, \alpha=0.3,0.6,0.9$ and $h=0.2,0.02$.
The figures of exact solution and the Caputo gH-derivative are shown
in Fig. \ref{f2}.


\begin{table}
\caption{Numerical results of Example \ref{ex2} for various $t,
\alpha$ and $h$.}\label{t2}
 \centering
\scalebox{0.6}{
\begin{tabular}{|c|l|l|l|l|l|l|}
  \hline
  \multicolumn{3}{ |c| }{$\alpha =0.3$}& \multicolumn{2}{ |c| }{$\alpha =0.6$ }& \multicolumn{2}{ |c| }{$\alpha =0.9$ } \\
 \hline
$t$   & $h=0.2$ & $h=0.02$ &  $h=0.2$ & $h=0.02$ &   $h=0.2$ & $h=0.02$ \\
\hline
0.1&(0,0.312475,0.624949)&(0,0.655421,1.31084)&(0,0.573896,1.14779)&(0,0.892967,1.78593)&(0,0.755737,1.51147)&(0,0.969249,1.9385)\\
0.2&(0,0.0976404,0.195281)&(0,0.429577,0.859154)&(0,0.329356,0.658712)&(0,0.797391,1.59478)&(0,0.571138,1.14228)&(0,0.939444,1.87889)\\
0.3&(0,0.0305101,0.0610203)&(0,0.281554,0.563107)&(0,0.189016,0.378032)&(0,0.712044,1.42409)&(0,0.43163,0.863261)&(0,0.910555,1.82111)\\
0.4&(0,0.00953365,0.0190673)&(0,0.184536,0.369072)&(0,0.108476,0.216951)&(0,0.635832,1.27166)&(0,0.326199,0.652398)&(0,0.882555,1.76511)\\
0.5&(0,0.00297902,0.00595805)&(0,0.120949,0.241898)&(0,0.0622536,0.124507)&(0,0.567777,1.13555)&(0,0.246521,0.493042)&(0,0.855415,1.71083)\\
0.6&(0,0.000930869,0.00186174)&(0,0.0792725,0.158545)&(0,0.0357271,0.0714542)&(0,0.507007,1.01401)&(0,0.186305,0.37261)&(0,0.829111,1.65822)\\
0.7&(0,0.000290873,0.000581746)&(0,0.0519568,0.103914)&(0,0.0205036,0.0410072)&(0,0.45274,0.905481)&(0,0.140797,0.281595)&(0,0.803615,1.60723)\\
0.8&(0,0.0000908904,0.000181781)&(0,0.0340536,0.0681072)&(0,0.0117669,0.0235339)&(0,0.404282,0.808565)&(0,0.106406,0.212812)&(0,0.778903,1.55781)\\
0.9&(0,0.000028401,0.0000568019)&(0,0.0223195,0.0446389)&(0,0.00675299,0.013506)&(0,0.361011,0.722022)&(0,0.0804149,0.16083)&(0,0.754951,1.5099)\\
1.0&(0,$8.87458\times10^{-6}$,0.0000177492)&(0,0.0146286,0.0292573)&(0,0.00387551,0.00775103)&(0,0.322371,0.644742)&(0,0.0607725,0.121545)&(0,0.731735,1.46347)\\
 \hline
 \end{tabular}
 }
\end{table}
\end{example}
\begin{example}\label{ex3} Let us to consider the following problem
$$
^{FC}D^{\alpha}_{\ast}y(t) = -\frac{\pi^2 t^{2-\alpha}
\alpha^2}{(2-3\alpha+\alpha^2)\Gamma(1-\alpha)}  {}_pF_q\left(1;\left[\frac{3}{2}-\frac{\alpha}{2},2-\frac{\alpha}{2}\right];-\frac{1}{4}
\pi^2t^2\alpha^2\right)\odot\left(0,\frac{1}{2},1\right), ~~1\leq
t \leq 2,
$$
where $y(1)=\left(0,\frac{1}{2},1\right)\odot\cos (\alpha\pi)$ and
the exact solution is $y(t)=\left(0,\frac{1}{2},1\right)\odot\cos
(\alpha\pi t)$. We know that this problem has the switching point at
$t=1.40426$. According to Eq.  (\ref{eq:5.21}) we should divide the
interval $[1,2]$ to the $N$ subinterval $[t_k, t_{k+1}]$, for $k=0, 1, ...,N-1$, and assuming the switching point belongs to $[t_j , t_{j+1}]$,  then the following iterative
formulas are applied as
$$
\begin{array}{l}
\displaystyle  y_{k+1} = y_k \oplus
\frac{h^{\alpha}}{\Gamma(\alpha+1)}\left(-\frac{\pi^2 t^{2-\alpha}
\alpha^2}{(2-3\alpha+\alpha^2)\Gamma(1-\alpha)}
~_pF_q\left(a; b; z t_k^2\right)\odot\left(0,\frac{1}{2},1\right)\right),\\~~~~~~k=0,1,\cdots, j,  \\
 \displaystyle  y_{k+1} = y_k \ominus (-1)
\frac{h^{\alpha}}{\Gamma(\alpha+1)}\left(-\frac{\pi^2 t^{2-\alpha}
\alpha^2}{(2-3\alpha+\alpha^2)\Gamma(1-\alpha)}~_pF_q\left(a; b; z
t_k^2\right)
\odot\left(0,\frac{1}{2},1\right)\right),\\~~~~~~k=j+1,\cdots, N-1,
\end{array}
$$
where $a=1,
b=\left[\frac{3}{2}-\frac{\alpha}{2},2-\frac{\alpha}{2}\right]$ and
$z=-\frac{1}{4} \pi^2 \alpha^2$. Numerical results are demonstrated
in Table \ref{t3} for $\alpha=0.8$ and $h=0.2,0.02,0.002$. In Fig.
\ref{f3}, the graphs of the exact solution and the Caputo
gH-derivative are presented for $\alpha=0.8$.
\begin{table}
\caption{Numerical results of Example \ref{ex3} for $\alpha=0.8$,
$h=0.2,0.02,0.002$ and various $t$.}\label{t3}
 \centering
\scalebox{0.9}{
\begin{tabular}{|c|l|l|l|}
  \hline
  \multicolumn{4}{ |c| }{$\alpha =0.8$}\\
 \hline
$t$   & $h=0.2$ & $h=0.02$ &  $h=0.002$ \\
\hline
1.1&(0,-0.452376,-0.918699)&(0,-0.463549,-0.928996)&(0,-0.464888,-0.929776)\\
1.2&(0,-0.489654,-0.984721)&(0,-0.495423,-0.991934)&(0,-0.496057,-0.992115)\\
1.3&(0,-0.489654,-0.984721)&(0,-0.495423,-0.991934)&(0,-0.496057,-0.992115)\\
1.4&(0,-0.452376,-0.918699)&(0,-0.463549,-0.928996)&(0,-0.464888,-0.929776)\\
1.5&(0,-0.400112,-0.800241)&(0,-0.404397,-0.808832)&(0,-0.404508,-0.809017)\\
1.6&(0,-0.307754,-0.628932)&(0,-0.317689,-0.637143)&(0,-0.318712,-0.637424)\\
1.7&(0,-0.20878,-0.417784)&(0,-0.21233,-0.425584)&(0,-0.21289,-0.425779)\\
1.8&(0,-0.0927856,-0.176232)&(0,-0.0935876,-0.187196)&(0,-0.0936907,-0.187381)\\
1.9&(0,0.0304523,0.0618529)&(0,0.0313271,0.0627734)&(0,0.0313953,0.0627905)\\
2.0&(0,0.146488,0.301241)&(0,0.15359,0.308805)&(0,0.154508,0.309017)\\
 \hline
 \end{tabular}
 }
\end{table}
\end{example}
\begin{example}\label{ex4}
Consider the following nonlinear fractional differential equations
under uncertainty:
$$
\sqrt{\eta}\odot {}^{FC}D^{\alpha}_{\ast}y(t) + y^2(t) = g(x), ~0<\alpha<1, ~t \in
[0,1],
$$
where
$$
g(x)=\left[ \frac{\Gamma (6)}{\Gamma(6-\alpha)} t^{5-\alpha}-\frac{3
\Gamma(5)}{\Gamma(5-\alpha)}t^{3-\alpha}+\frac{
\Gamma(5)}{\Gamma(4-\alpha)}t^{3-\alpha}+(t^5-3t^4+2t^3)^2
\right]\odot\tilde{\eta},
$$
and $\tilde{\eta} (r) = (1,2,3),~y\neq 0$. Then the exact
solution of the problem is $y(t) = \sqrt{\eta}\odot
(t^5-3t^4+2t^3)$. By solving the problem under
$^{FC}$[(i)-gH]-differentiability using
$$
y_{k+1}=y_k+\frac{h^\alpha}{\Gamma(\alpha+1)}\odot \Gamma_k,~~~~
k=0,1,\cdots, N-1,
$$
we obtain the numerical solution shown in Table \ref{t4}, with
different order of differentiability and step size.
\begin{table}
\caption{Absolute error of Example \ref{ex4} at $t=1$.}\label{t4}
 \centering
\scalebox{0.9}{
\begin{tabular}{|c|l|l|l|l|l|}
  \hline
$h$   & $\alpha=0.1$ & $\alpha=0.3$ &  $\alpha=0.5$ &  $\alpha=0.7$ &  $\alpha=0.9$\\
\hline
$\frac{1}{10}$&$7.20602\times 10^{-2}$&$6.5418\times 10^{-2}$&$5.8823\times 10^{-2}$&$5.3707\times 10^{-2}$&$5.0201\times 10^{-2}$\\
$\frac{1}{20}$&$3.9603\times 10^{-2}$&$3.3498\times 10^{-2}$&$2.9368\times 10^{-2}$&$2.6937\times 10^{-2}$&$2.5668\times 10^{-2}$\\
$\frac{1}{40}$&$2.0653\times 10^{-2}$&$1.6611\times 10^{-2}$&$1.4420\times 10^{-2}$&$1.3384\times 10^{-2}$&$1.2962\times 10^{-2}$\\
$\frac{1}{80}$&$1.0448\times 10^{-4}$&$8.1038\times 10^{-3}$&$7.0482\times 10^{-3}$&$6.6381\times 10^{-3}$&$6.5091\times 10^{-3}$\\
 \hline
 \end{tabular}
 }
\end{table}
 In Fig.\ref{f4}, the graphs of the exact solution and Caputo
gH-derivatives are presented for $\alpha=0.1, 0.3, 0.5, 0.7, 0.9, 1$ and in Fig.\ref{f5} these Caputo
gH-derivatives have been compared in $r=0.5$.



\end{example}
\begin{remark}Although, we have obtained the solution under $^{FC}[(i)-gH]$-differentiability, but it is easy to check that it is not $^{FC}[(i)-gH]$-differentiable on $(0,1)$. Actually, due to obtained results ( see Table \ref{t5}), we can consider the proper interval that the given exact solution and its approximation is $^{FC}[(i)-gH]$-differentiable. Also, note that, we have computed the approximation of the solution of Example \ref{ex4} at point $t=1$, which is clearly this point take place out of proper domain of $^{FC}[(i)-gH]$-differentiability. In fact, the computed error at point $t=1$, just obtained based on the lower-upper approximation of lower-upper of exact solution. For more clarification, we determined switching points regarding each order of differentiability.
\end{remark}
\begin{table}
\caption{Switching points for different values of $\alpha$}\label{t5}
 \centering
\scalebox{0.9}{
\begin{tabular}{|c|l|l|l|l|l|l|}
  \hline
$\alpha$& $0.1$& $0.3$& $0.5$& $0.7$& $0.9$& $1$\\
\hline
$t$&$0.9701$&$0.9109$&$0.8525$&$0.7949$&$0.7381$&$0.7101$\\
 \hline
 \end{tabular}
 }
\end{table}
\begin{remark}
Indeed, using results of Table \ref{t5}, in fact, we deduce that by considering the problem of fractional order instead of integer order (here, first order), we obtain some wider interval Than the first order case, on the other hand, when $\alpha=1$, the valid interval that the given exact solution verify the assumption $^{FC}[(i)-gH]$-differentiability is $(0,0.7101)$, while for $\alpha=0.9$ and $\alpha=0.7$, the valid interval are $(0, 0.7381)$ and $(0, 0.7949)$, respectively. Actually, this is a first time in the literature that this new results, i.e., extending the length of valid interval that the type of differentiability remains unchanged, is investigated.
\end{remark}
\section{Conclusions}
Fractional differential equations are one of the important topics of
the fuzzy arithmetic which have many applications in sciences and
engineering. Thus finding the numerical and analytical methods to
solve these problems is very important. This paper was presented
based on the two main topics. Firstly, proving the generalized
Taylor series expansion for fuzzy valued function based on the
concept of generalized Hukuhara differentiability. Secondly,
introducing the fuzzy generalized Euler's method as an applications
of the generalized Taylor expansion and applying it to solve the
fuzzy fractional differential equations. The capabilities and
abilities of presented method were showed by presenting several
theorem about the consistence, the convergence and the stability of
the generalized Euler's method. Also, accuracy and efficiency of
method were illustrated by considering on the local and global
truncation errors. The numerical results specially in the switching
point case showed the precision of the generalized Euler's method to
solve the fuzzy fractional differential equations.
\section*{Funding}
The work of J. J. Nieto has been partially supported by Agencia
Estatal de Investigaci\'on (AEI) of Spain under grant
MTM2016-75140-P,
 co-financed by the European Community fund FEDER, and XUNTA de Galicia under grants GRC2015-004 and R2016-022.
\section*{Acknowledgements}
The authors are grateful to the anonymous reviewers for their
helpful, valuable comments and suggestions in the improvement of
this manuscript.

\end{document}